\documentclass[12pt]{elsarticle}

\usepackage{amsmath,amssymb}
\usepackage{mathrsfs}
\usepackage{makecell}
\usepackage{url}
\usepackage[english]{babel}
\usepackage{wrapfig}
\usepackage{color}
\usepackage{array}
\usepackage{multirow}
\usepackage{bigstrut}
\usepackage{booktabs}
\usepackage{txfonts}
\usepackage{natbib}
\usepackage{bm}

\newcommand{\ii}{\textnormal{i}}
\newcommand{\ee}{\textnormal{e}}

\newcommand{\len}{4cm}
\newcommand{\lenx}{-2.4cm}
\newcommand{\lem}{4cm}
\newcommand{\lemx}{-2.4cm}

\newcommand{\eee}{\mathbf{e}}
\newcommand{\fff}{\mathbf{f}}
\newcommand{\xx}{\bm{x}}
\newcommand{\qq}{\bm{q}}
\newcommand{\kk}{\bm{k}}

\newcommand{\uu}{\bm{u}}
\newcommand{\eps}{\mbox{\boldmath$\varepsilon$}}
\newcommand{\sig}{\mbox{\boldmath$\sigma$}}

\newcommand{\tautau}{\mbox{\boldmath$\tau$}}

\newcommand{\real}{\operatorname{Re}}
\newcommand{\imag}{\operatorname{Im}}

\newcommand{\mysym}[1]{{#1}_\textnormal{sym}\hspace{-0.2cm}}

\def\ack{\section*{Acknowledgements}
  \addtocontents{toc}{\protect\vspace{6pt}}
  \addcontentsline{toc}{section}{Acknowledgements}
}

\journal{CRAS}

\begin{document}

\begin{frontmatter}

\title{Fourier-based schemes for computing the mechanical response of composites with accurate local fields}

\author{Fran\c{c}ois Willot\footnote{Tel.: +33 1 64 69 48 07.}}
\ead{francois.willot@ensmp.fr}
\address{Mines ParisTech, PSL Research University, Centre for Mathematical Morphology,\\
35, rue St Honor\'e, 77300 Fontainebleau, France}

\begin{abstract}
We modify the Green operator involved in Fourier-based computational schemes in 
elasticity, in 2D and 3D.
The new operator is derived by expressing continuum mechanics in terms of centered differences
on a rotated grid.
Use of the modified Green operator leads, in all systems investigated, to more accurate strain and stress fields
than using the discretizations proposed by Moulinec and Suquet~(1994) or Willot and Pellegrini~(2008).
Moreover, we compared the convergence rates of the ``direct'' and ``accelerated'' FFT schemes with the different discretizations.
The discretization method proposed in this work allows for much faster 
FFT schemes with respect to two criteria: stress equilibrium and effective elastic moduli.
\end{abstract}


\end{frontmatter}

\section{Introduction}\label{sec:intro}
Fourier-based algorithms, or ``FFT'' methods for short, are an efficient approach
for computing the mechanical response of composites.
Initially restricted to linear-elastic media,
FFT tools are nowadays employed to treat more involved problems, ranging from 
viscoplasticity~\cite{lee2011modeling} to crack propagation~\cite{li2011non}.
In FFT methods, 
the microstructure is defined by 2D or 3D images 
and the local stress and strain tensors are computed along each pixel or ``voxel'' in the image. 
Coupled with automatic or semi-automatic image segmentation techniques~\cite{faessel13}, this allows
for the computation of 
the mechanical response of a material
directly from experimental acquisitions, like focused ion beam 
or 3D microtomography techniques~\cite{Willot13a}.
The latter often deliver images containing billions of voxels, for which
FFT methods are efficient~\cite{Dunant13,Willot10}.
This allows one to take into account representative volume elements 
of materials which are multiscale by nature such as concrete or mortar~\cite{Escoda11}.
From a practical viewpoint, the simplicity of FFT methods is attractive to
researchers and engineers who need not be experts in
the underlying numerical methods to use them.
Nowadays, FFT tools are available not only as academic or free
softwares~\cite{morphhom,craft} but also as commercial ones~\cite{geodict}.

In the past years, progress has been made in the understanding of FFT algorithms.
Vond\v{r}ejc and co-workers have recently shown that the original method of
Moulinec and Suquet~\cite{Moulinec94xx}
corresponds, under 
one technical assumption, to a particular choice of approximation space and optimization method~\cite{zeman14}
(see also~\cite{Brisard12}).
This property allows one to derive other FFT schemes that use standard optimization algorithms,
such as the conjugate gradient method.
In this regard, making use of variational formulations, 
efficient numerical methods have been proposed
that combine FFTs with an underlying gradient descent algorithm~\cite{Zeman10,Brisard10}.

Different approximation spaces or discretization methods have also been proposed, where,
contrary to the original scheme, the fields are not trigonometric polynomials anymore.
Brisard and Dormieux introduced ``energy-based'' FFT schemes 
that rely on Galerkin approximations of Hashin and Shtrikman's
variational principle~\cite{Brisard10,Brisard12} and derived
modified Green operators
consistent with the new formulation.
They obtained improved convergence properties and 
local fields devoid of the spurious oscillations
observed in the original FFT scheme~\cite{Brisard12,Brisard14}.
In the context of conductivity,
accurate local fields and improved convergence rates
have also been obtained
from modified Green operators based on finite-differences~\cite{Willot14a}.
These results follow
earlier works where continuum mechanics
are expressed by centered~\cite{Muller96,Brown02} or ``forward and backward'' finite differences~\cite{willot08c}.

This work focuses on the effect of discretization in FFT methods.
It is organized as follows.
We first recall the equations of elasticity in the continuum
(Sec.~\ref{sec:problem}).
We give the Lippmann-Schwinger equations and the ``direct'' and ``accelerated''
FFT schemes in Sec.~(\ref{sec:greenCont}).
In Sec.~(\ref{sec:disc}), a general formulation of the Green operator is derived
that incorporates methods in~\cite{willot08c}, and
a new discretization scheme is proposed.
The accuracy of the local stress and strain fields are examined in Sec.~(\ref{sec:accuracy})
whereas the convergence rates of the various FFT methods are investigated in Sec.~(\ref{sec:rate}).
We conclude in Sec.~(\ref{sec:conclusion}).

\section{Microstructure and material elastic response}\label{sec:problem}
We are concerned with solving the equations of linear elasticity 
in a square or cubic domain $\Omega=[-1/2;1/2]^d$ in dimension $d$ ($d=2$ or $3$):
\begin{equation}\label{eq:a1}
\sigma_{ij}(\xx)=C_{ij,kl}(\xx)\varepsilon(\xx), \quad 
\partial_i \sigma_{ij}(\xx)\equiv 0, \quad
\varepsilon_{ij}(\xx)=(1/2)\left[\partial_i u_j(\xx)+\partial_j u_i(\xx)\right], \quad
\end{equation}
where $\eps(\xx)$ is the strain field,
$\sig(\xx)$ the stress field,
$\uu(\xx)$ the displacement vector field,
$\mathbb{C}(\xx)$ the local elasticity tensor
and $\xx$ is a point in $\Omega$.
Tensorial components refer to a system of Cartesian coordinates $(\eee_1;\eee_2)$ in 2D
and $(\eee_1;\eee_2;\eee_3)$ in 3D.
The material has an isotropic local elastic response that reads:
\begin{equation}
C_{ij,kl}(\xx)=\lambda(\xx)\delta_{ij}\delta_{kl}
+\mu(\xx)(\delta_{ik}\delta_{jl}+\delta_{il}\delta_{jk}),
\end{equation}
where $\delta$ is the Kronecker symbol
and $\lambda(\xx)$ and $\mu(\xx)$ are constant-per-phase Lam\'e's first and second coefficients.
The local bulk modulus $\kappa=\lambda+(2/d)\mu$ and
the elastic moduli take on
values:
$$
  \lambda(\xx)=\lambda^\alpha,
  \quad \kappa(\xx)=\kappa^\alpha,
  \quad \mu(\xx)=\mu^\alpha,
$$
in phase $\alpha$.
For simplicity, we restrict ourselves to binary media and, by convention,
$\alpha=1$ is the matrix and $\alpha=2$ the inclusions.
Hereafter, we fix Poisson's ratios in each phase to $\nu^1=\nu^2=0.25$
so that, in 3D and 2D~\cite{thorpe84}, we have $\mu^\alpha/\kappa^\alpha=0.6$.
The contrast of properties $\chi$ reads:
\begin{equation}
\chi=\frac{\kappa^2}{\kappa^1}=\frac{\mu^2}{\mu^1}=\frac{\lambda^2}{\lambda^1},
\end{equation}
where $0\leq\chi\leq\infty$.
In the matrix, we also fix $\kappa^1=1$ ($d=2$ or $3$), $\mu^1=0.6$ ($d=2$ or $3$),
 $\lambda^1=0.4$  ($d=2$), $\lambda^1=0.6$  ($d=3$), so that
the local properties of the material are parametrized by one unique variable,
the contrast of properties $\chi$.
In 3D, the Young modulus is $E^1=3/2$ in the matrix
and $E^2=3\chi/2$ in the inclusion.
The medium is porous when $\chi=0$ and rigidly-reinforced
when $\chi=\infty$.

Periodic boundary conditions are applied with the material 
subjected to an overall strain loading $\overline{\eps}$:
\begin{equation}\label{eq:a2}
\sigma_{ij}(\xx) n_j(\xx) -\# \quad (\xx\in\partial\Omega),\qquad
\langle\varepsilon_{kl}(\xx)\rangle=\overline{\varepsilon_{kl}},
\end{equation}
where $\bm{n}$ is the normal at the boundary $\partial\Omega$ of the domain $\Omega$, oriented outward,
$-\#$ denotes anti-periodicity and $\langle\cdot\rangle$ denotes the spatial mean over $\Omega$.
The resulting effective elastic tensor $\widetilde{\mathbb{C}}$ is computed from:
\begin{equation}
\langle\sigma_{ij}(\xx)\rangle= \widetilde{C}_{ij,kl}\overline{\varepsilon_{kl}}.
\end{equation}

\section{Lippmann-Schwinger equation and FFT methods}\label{sec:greenCont}
Fourier methods are by principle based on the Lippmann-Schwinger equations.
The latter follow from~(\ref{eq:a1}) and~(\ref{eq:a2}) as~\cite{MiltonBook}:
\begin{equation}\label{eq:ls}
\tau_{ij}(\xx)=\sigma_{ij}(\xx)-C^0_{ij,kl}\varepsilon_{kl}(\xx), \quad
\varepsilon_{ij}(\xx)=\overline{\varepsilon_{ij}}-
 \int_{\xx'}{\rm d}^d\xx'\,G_{ij,kl}(\xx'-\xx)\tau_{kl}(\xx'),
\end{equation}
where we have introduced a homogeneous ``reference'' elasticity tensor $\mathbb{C}^0$
and its associated polarization field $\tautau$
and Green operator $\mathbb{G}$.
In the above we assume $\langle\mathbb{G}\rangle=0$ so that 
$\overline{\eps}=\langle\eps(\xx)\rangle$ holds.
The Green operator has, in the Fourier domain, the closed form~\cite{KanaunBook}:
\begin{equation}\label{eq:green0}
G_{ij,kl}(\qq)=\left\lbrace
 q_i\left[q_mC^0_{mj,kn} q_n\right]^{-1} q_l
 \right\rbrace_\textnormal{sym},
\end{equation}
where $\qq\neq 0$ are the Fourier wave vectors and the subscript $_\textnormal{sym}$ indicates minor symmetrization
with respect to the variables $(i,j)$ and $(k,l)$.
Hereafter, we assume that $\mathbb{C}^0$ is 
a symmetric, positive-definite, isotropic tensor defined by its bulk ($\kappa^0$) and shear ($\mu^0$) moduli,
 or Lam\'e coefficient ($\lambda^0$).
Accordingly, when $\qq\neq 0$, the Green operator $\mathbb{G}$ is also symmetric definite and
we have:
\begin{equation}\label{eq:green0b}
G_{ij,kl}(\qq)=
\frac{1}{\mu^0}\left[
  \left(\frac{q_iq_l}{|\qq|^2}\delta_{jk}\right)_\textnormal{sym}
  -\frac{\lambda^0+\mu^0}{\lambda^0+2\mu^0}\frac{q_i q_j q_k q_l}{|\qq|^4}
   \right].
\end{equation}

The ``direct scheme''~\cite{Moulinec94xx} consists in applying Eqs.~(\ref{eq:ls}) iteratively as:
\begin{equation}\label{eq:ds}
  \eps^{k=0}\equiv\overline{\eps}, \qquad
  \eps^{k+1}=\overline{\eps}-\mathbb{G}\ast (\sig-\mathbb{C}^0:\eps^k) \quad (k\geq 0).
\end{equation}
In Moulinec and Suquet's method, the convolution product ($\ast$) above is computed as an algebraic product in the Fourier domain,
making use of~(\ref{eq:green0b}).
Discrete Fourier transforms are used to switch between the space $\Omega$ and Fourier domain $\mathcal{F}$.
This amounts to representing the strain field as a trigonometric polynomial~\cite{zeman14} of the form:
\begin{equation}\label{eq:ft}
\eps(\xx)=\frac{1}{L^d}\sum_{\qq\in\mathcal{F}} \eps(\qq)\,\ee^{\ii \qq\cdot\xx},
\end{equation}
where $\eps(\qq)$ is the discrete Fourier transform of $\eps(\xx)$. Accordingly:
\begin{equation}\label{eq:ft2}
\eps(\qq)=\sum_{\xx\in\Omega} \eps(\xx)\,\ee^{-\ii \qq\cdot\xx}.
\end{equation}
Similar forms are used for the stress and displacement fields.
In practice, the domain $\Omega$ is discretized on a square or cubic grid of $L^d$ voxels and
the operator $\mathbb{G}(\qq)$ in~(\ref{eq:green0}) is evaluated 
along equispaced Fourier modes:
\begin{equation}\label{eq:contmodes}
q_i=\frac{2\pi m_i}{L},\qquad  m_i=
\left\lbrace\begin{array}{cc}
 1-(L/2), ..., L/2,   & (L\textnormal{ even}),\\
 -(L-1)/2, ..., (L-1)/2, & (L\textnormal{ odd}).
 \end{array}\right.
\end{equation}
As noted in~\cite{moulinec1998numerical}, when $L$ is even,
the relation:
\begin{equation}
\mathbb{G}(\qq)^*=\mathbb{G}(-\qq),
\end{equation}
where $\mathbb{G}(\qq)^*$ is the complex conjugate of $\mathbb{G}(\qq)$,
is not verified when one of the components of $\qq$ is equal to the highest frequency
$q_i=\pi$ (i.e. $m_i=L/2$).
As a consequence, the backward Fourier transform
of $\mathbb{G}(\qq)\tautau(\qq)$ used to compute the strain field has non-zero imaginary part even if 
$\tautau(\xx)$ is purely real.
To fix this problem, we follow~\cite{moulinec1998numerical} and set:
\begin{equation}\label{eq:hfreq}
 \mathbb{G}(\qq)=\left(\mathbb{C}^0\right)^{-1}, \qquad \textnormal{ if $q_i=\pi$ for some $i$.}
\end{equation}
This choice enforces $\sig(\qq)=0$ at the concerned Fourier modes.
In doing so, the strain field $\eps$ in~(\ref{eq:ft}) is not strictly-speaking irrotational, because
of the lack of constraint at high Fourier modes for the strain field.
We briefly mention another option that we explored in this work.
It consists in forcing the symmetry by replacing $\mathbb{G}(\qq)$ with:
\begin{equation}\label{eq:forcing}
  \frac{\mathbb{G}(\qq)^*+\mathbb{G}(-\qq)}{2},
\end{equation}
when one of the components of $\qq$ equals $\pi$, which enforces $\eps(\qq)=0$ at the highest modes.
Numerical experiments indicate that the choice for $\mathbb{G}(\qq)$
at the highest frequencies has little influence on the convergence rate,
except at small resolution. When $L<128$ pixels, faster convergence was achieved with the choice
$\mathbb{G}(\qq)=\left(\mathbb{C}^0\right)^{-1}$.
Furthermore, in the 2D example studied in this work,
the choice $\mathbb{G}(\qq)=\left(\mathbb{C}^0\right)^{-1}$
led to smaller oscillations, consistently with observations in~\cite{moulinec1998numerical}.
The use of~(\ref{eq:forcing}) was therefore not pursued further.
We emphasize that, when $L$ is odd, this discrepancy disappears and no special treatment is needed.

Refined FFT algorithms have been introduced to overcome the slow convergence rate
of the direct scheme, observed for highly-contrasted composites,
most notably the ``accelerated scheme''~\cite{Eyre99} and ``augmented Lagrangian''~\cite{Michel01}
methods.
In this work, we use the extension of the accelerated scheme
to elasticity~\cite{Michel01,moulinec13}:
\begin{equation}\label{eq:as}
 \eps^{k+1}=\eps^k+2\left( \mathbb{C}+\mathbb{C}^0 \right)^{-1}:\mathbb{C}^0:\left[
 \overline{\eps}-\eps^k-\mathbb{G}\ast\left(\mathbb{C}:\eps^k-\mathbb{C}^0:\eps^k\right)
 \right].
\end{equation}
The convergence rates of the accelerated and direct schemes
depend on the choice of the reference tensor $\mathbb{C}^0$.
For the accelerated scheme the optimal choice is~\cite{Eyre99,moulinec13}:
\begin{equation}\label{eq:refAcc}
\kappa^0= \sqrt{\kappa^1\kappa^2},\qquad
\mu^0= \sqrt{\mu^1\mu^2}.
\end{equation}
For the direct scheme, upper bounds on the eigenvalues of the Green operator
suggest the choice~\cite{moulinec13}:
\begin{equation}\label{eq:ref}
\kappa^0=\beta\left(\kappa^1+\kappa^2\right),\qquad
\mu^0=\beta\left(\mu^1+\mu^2\right), 
\end{equation}
with $\beta=1/2$. 

\section{Discretization and approximation space}\label{sec:disc}
In this section, we derive the expression of a modified Green operator $\mathbb{G}'$
that replaces $\mathbb{G}$
defined in~(\ref{eq:green0b}) and~(\ref{eq:contmodes}).
We give it in a form that includes previously-proposed modified operators~\cite{willot08c}
and also introduce a new one.

\subsection{Two dimensions}
In the following, we assume that the strain and stress fields are defined
on a grid of points in 2D, one per pixel.
Eqs.~(\ref{eq:ft}) and~(\ref{eq:ft2}) are used to apply discrete Fourier transforms,
but we do not postulate a representation in the continuum anymore.
The equilibrium and strain admissibility conditions~(\ref{eq:a1})
are approximated by means of finite differences on which we apply the discrete transforms~(\ref{eq:ft}) and~(\ref{eq:ft2}).
In~\cite{willot08c},
this results in the following form:
\begin{equation}\label{eq:discFour}
k_i^*(\qq)\sigma_{ij}(\qq)=0, \qquad 
\varepsilon_{ij}(\qq)=(1/2)\left[ k_i(\qq)u_j(\qq) + k_j(\qq)u_i(\qq) \right],
\end{equation}
where $\kk$ and $\kk^*$ represent ``discrete'' gradient and divergence operators, respectively.
In the centered scheme, one takes $\kk$ equal to:
\begin{equation}\label{eq:mul2}
k^\textnormal{C}_i(\qq)=\ii \sin\left(q_i\right),
\end{equation}
whereas in scheme~\cite{willot08c}, one chooses, for $\kk$:
\begin{equation}\label{eq:wil072}
k^\textnormal{W}_i(\qq)=\ee^{\ii q_i}-1.
\end{equation}
These expressions correspond, respectively, to the centered scheme:
\begin{equation}\label{eq:mul1}
\partial_j\sigma_{ij}(\xx)\approx \frac{\sigma_{ij}(\xx+\ee_j)-\sigma_{ij}(\xx-\ee_j)}{2}, \qquad
\partial_j u_i(\xx)\approx \frac{u_i(\xx+\ee_j)-u_i(\xx-\ee_j)}{2},
\end{equation}
and to the forward-and-backward difference scheme:
\begin{equation}\label{eq:wil071}
\partial_j\sigma_{ij}(\xx)\approx \sigma_{ij}(\xx)-\sigma_{ij}(\xx-\ee_j), \qquad
\partial_j u_i(\xx)\approx u_i(\xx+\ee_j)-u_i(\xx).
\end{equation}
Using~(\ref{eq:discFour}), the resulting Green operator reads:
\begin{equation}\label{eq:discFour2}
G'_{ij,kl}(\qq)=\left\lbrace
 k_i(\qq) \left[k_m(\qq) C^0_{mj,kn} k_n^*(\qq)\right]^{-1} k_l^*(\qq)
 \right\rbrace_\textnormal{sym},
\end{equation}
which is homogeneous in $\kk$, so that, with $r_i=k_i/|\kk|$:
\begin{equation}\label{eq:green1}
G'_{ij,kl}(\qq)=
\frac{
   (\lambda^0+2\mu^0)\left(r_ir_l^*\delta_{jk}\right)_\textnormal{sym}
        -\lambda^0\real\left(r_i r_j^*\right)\real\left(r_k r_l^*\right)-\mu^0 r_i r_j\left(r_k r_l\right)^*
 }{
   \mu^0\left[2(\lambda^0+\mu^0)-\lambda^0\left|r_1^2+r_2^2\right|^2\right]
 },
\end{equation}
where $\real(\cdot)$ denotes the real part of the enclosed quantity.
The denominator on the right-hand side is strictly positive due to the triangle inequality.
Hereafter, the operator $\mathbb{G}'$ is denoted 
by $\mathbb{G}^\textnormal{C}$ when $\kk=\kk^\textnormal{C}$ 
and by $\mathbb{G}^\textnormal{W}$ 
when $\kk=\kk^\textnormal{W}$. 
Note that the operator $\mathbb{G}$ in~(\ref{eq:green0b})
is recovered from~(\ref{eq:green1}) by setting $k_i(\qq)=\ii q_i$.
Now the Green operator $\mathbb{G}'$ is complex and follows the minor and major symmetries:
\begin{equation}\label{eq:greenSym}
G'_{ij,kl}(\qq)=G'_{ji,kl}(\qq)=G'_{ij,lk}(\qq), \qquad
G'_{ij,kl}(\qq)=\left[G'_{kl,ij}(\qq)\right]^*.
\end{equation}
Furthermore we have:
\begin{equation}\label{eq:greenSym2}
\mathbb{G}^\textnormal{C,W}(\qq)^*=\mathbb{G}^\textnormal{C,W}(-\qq),
\end{equation}
for the schemes~(\ref{eq:mul1}) and~(\ref{eq:wil071}), including when $q_i=\pi$ and when $L$ is even,
because $k_i^{C,W}$ is real at the frequency $q_i=\pi$. Therefore,
the fix~(\ref{eq:hfreq}) in Sec.~(\ref{sec:greenCont}) is not necessary.
However, a problem of a different nature arises when using the centered scheme~(\ref{eq:mul2}) when $L$ is even.
Eq.~(\ref{eq:green1}) does not define the Green operator $\mathbb{G}^\textnormal{C}$
at the three frequencies $\qq=(\pi;0)$, $(0;\pi)$ and $(\pi;\pi)$, for which $\kk^\textnormal{C}(\qq)=0$.
This is because the second equation in~(\ref{eq:mul1}) has in general non-unique
solutions for the displacement field $\uu(\xx)$.
Indeed, when $L$ is even, the displacement is defined up to a linear combination
of $2$-voxels periodic fields.
They are given by the following $8$ independent fields:
\begin{eqnarray}\label{eq:mulDisp}
&
v^1_m(\xx)=\delta_{mn},\quad
v^2_m(\xx)=(-1)^i\delta_{mn},\quad
v^3_m(\xx)=(-1)^j\delta_{mn},\quad
v^4_m(\xx)=(-1)^{i+j}\delta_{mn},\quad
\nonumber &\\ & \hspace{5cm}
\xx=\left(\frac{i}{L};
          \frac{j}{L}\right),\,\,\,\,
i, j=1-\frac{L}{2}, ..., \frac{L}{2},\,\,\,\,
n=1, 2. &
\end{eqnarray}
The operator $\mathbb{G}^\textnormal{C}$ remains finite 
when $\qq$ approaches one of the modes $(\pi;0)$,
$(0;\pi)$ or $(\pi;\pi)$, but can not be continuously extended at these modes.
To fix this problem, we set, for the centered scheme:
\begin{equation}\label{eq:mulfix}
G^\textnormal{C}_{ij,kl}(\qq)=0, \qquad
  \textnormal{ if $L$ is even and ($q_i=0$ or $q_i=\pi)$ for all $i$},
\end{equation}
which enforces $\eps(\qq)=0$ at the highest frequencies.
The strain field $\eps(\xx)$ is accordingly admissible,
and the stress field $\sig(\xx)$ is divergence-free, in the sense of~(\ref{eq:mul1}).
We explored the alternate choice $\mathbb{G}^\textnormal{C}(\qq)=\left(\mathbb{C}^0\right)^{-1}$
in~(\ref{eq:mulfix}).
Almost identical convergence rates and oscillations were observed for the two options,
however the choice $\mathbb{G}^\textnormal{C}(\qq)=\left(\mathbb{C}^0\right)^{-1}$
does not produce an irrotational strain field and is not considered further.
We emphasize that no special treatment is required
for the operator $\mathbb{G}^\textnormal{W}$ at high frequencies since $\kk^\textnormal{W}\neq 0$
when $\qq\neq 0$.

By substituting $\mathbb{G}$ with $\mathbb{G}^\textnormal{C}$
or $\mathbb{G}^\textnormal{W}$ in (\ref{eq:ds}) and (\ref{eq:as}),
we derive ``direct'' and ``accelerated'' schemes that solve~(\ref{eq:mul1}) or~(\ref{eq:wil071}).
In the limit of very fine resolution,
we have $\mathbb{G}^\textnormal{C,W}(\qq)\approx \mathbb{G}(\qq)$
when $\qq\to 0$, which
guarantees that the strain and stress fields do not depend on the employed discretization.
This property holds for any choice of $\kk$ such as $\kk\sim \ii\qq$ when $\qq\to 0$.

On the one hand, derivatives are estimated 
more locally in the forward-and-backward scheme~(\ref{eq:wil072}) than
in the centered scheme~(\ref{eq:mul2}), which is important along interfaces. On the other hand, 
the forward-and-backward scheme does not treat symmetrically the two angle bisectors $\eee_1+\eee_2$ and $\eee_1-\eee_2$~\cite{willot08c}.
In a domain containing a single centered disc, the scheme produces fields
that break the axial symmetries of the problem.
In fact, the discretization~(\ref{eq:wil071}) is actually one of four possible choices, all of them breaking the symmetries.
Attempts to force the symmetry by averaging over the four Green operators or over the fields themselves,
as proposed in~\cite{willot08c},
are not explored in this work.
The fomer method indeed leads to less accurate ``diffuse'' local fields.
The latter necessitates to run four
different computations, in 2D, instead of one, which is cumbersome.

\begin{figure}
\begin{center}
\includegraphics[width=4cm]{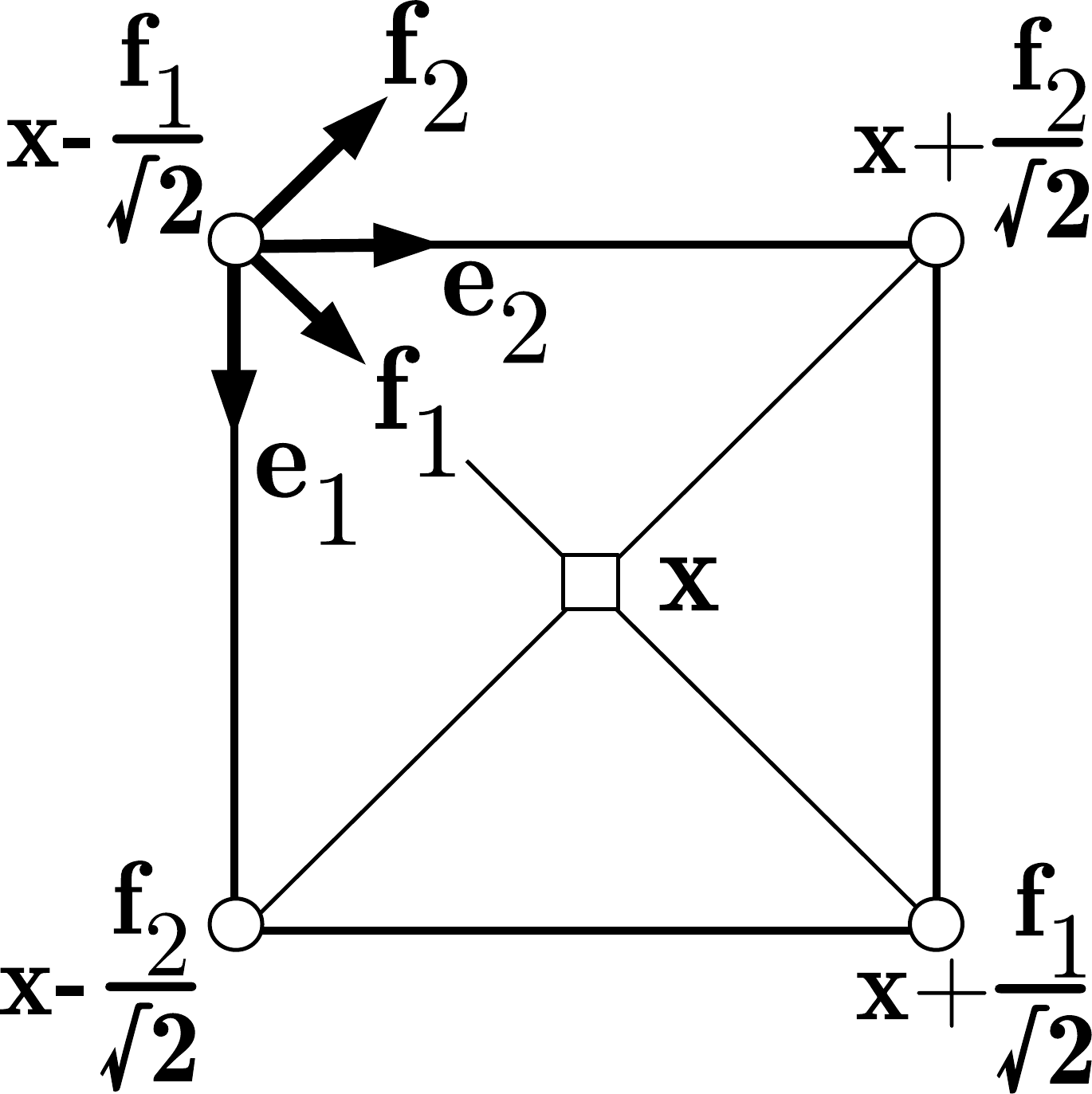}
\caption{\label{fig:2drot}
A pixel with edges parallel to the Cartesian axis $(\eee_1;\eee_2)$.
Superimposed: $45^\textnormal{o}$-rotated basis ($\fff_1; \fff_2)$.
The strain and stress fields are evaluated at the pixel center $\xx$ (square).
The displacement and the divergence of the stress field lie along the pixel corners (disks).
}
\end{center}
\end{figure}

In the rest of this section, we derive a discrete scheme in 2D different from~(\ref{eq:mul1}) and~(\ref{eq:wil071}).
In this scheme, the displacement field is evaluated at the $4$ corners of the pixels
and the strain and stress fields are evaluated at the centers of the pixels.
We first express these fields in the $45^\textnormal{o}$-rotated
basis:
\begin{equation}
  \fff_1=\frac{\eee_1+\eee_2}{\sqrt{2}}, \qquad \fff_2=\frac{\eee_2-\eee_1}{\sqrt{2}},
\end{equation}
by:
\begin{equation}
u_i=R_{iI}u_I,\quad \varepsilon_{ij}=R_{iI}\varepsilon_{IJ} R'_{Jj}, \quad \sigma_{ij}=R_{iI}\sigma_{IJ} R'_{Jj},
\qquad R_{iJ}=\frac{1-2\delta_{i1}\delta_{J2}}{\sqrt{2}},
\end{equation}
where uppercase indices refer to components in the rotated grid.
We discretize~(\ref{eq:a1}) in the rotated basis by the centered differences (see Fig.~\ref{fig:2drot}): 
\begin{subequations}\label{eq:discFull}
\begin{eqnarray}
&\sigma_{IJ}\left(\xx\right)= C_{IJ,KL}\left(\xx\right)
   \varepsilon_{KL}\left(\xx\right),&\\ \label{eq:discFullStrainComp}
&\sigma_{I1}\left(\xx\right)-\sigma_{I1}\left(\xx-\sqrt{2}\fff_1\right)
 +\sigma_{I2}\left(\xx+\frac{\fff_2-\fff_1}{\sqrt{2}}\right)-\sigma_{I2}\left(\xx-\frac{\fff_1+\fff_2}{\sqrt{2}}\right)=0, 
  &\\  \label{eq:discFullStressDiv}
&\varepsilon_{KL}\left(\xx\right)=\frac{1}{2\sqrt{2}}\left[
     u_K\left(\xx+\frac{\fff_L}{\sqrt{2}}\right)-u_K\left(\xx-\frac{\fff_L}{\sqrt{2}}\right) 
   + u_L\left(\xx+\frac{\fff_K}{\sqrt{2}}\right)-u_L\left(\xx-\frac{\fff_K}{\sqrt{2}}\right)\right]\label{eq:discFulleps}&
\end{eqnarray}
\end{subequations}
where $\xx$ lie at the centers of the pixels and 
$\xx\pm\fff_I/\sqrt{2}$ lie at the corners.
Expressing back~(\ref{eq:discFull}) in the original Cartesian grid $(\eee_1;\eee_2)$ 
and applying the backward discrete Fourier transform~(\ref{eq:ft}) we arrive 
again at~(\ref{eq:discFour}) with
the following expression for $\kk$: 
\begin{equation}\label{eq:wil14}
k_i^\textnormal{R}(\qq)=\frac{\ii}{2}\tan\left(\frac{q_i}{2}\right)\left(1+\ee^{\ii q_1}\right)\left(1+\ee^{\ii q_2}\right).
\end{equation}
We denote by $\mathbb{G}^\textnormal{R}$ the corresponding Green operator, derived by substituting $\kk=\kk^\textnormal{R}$
in~(\ref{eq:green1}).
The operator $\mathbb{G}^\textnormal{R}$ is real and also verifies:
\begin{equation}\label{eq:grsym}
\mathbb{G}^\textnormal{R}(\qq)=\mathbb{G}^\textnormal{R}(-\qq).
\end{equation}
However, when $L$ is even, 
$\kk^\textnormal{R}=0$ when $\qq=(\pi;\pi)$ and 
$\mathbb{G}^\textnormal{R}$ is not defined by~(\ref{eq:green1}) at this frequency.
Again, this is because (\ref{eq:discFulleps}) gives the displacement field
up to linear combinations of the $4$ independent fields $v^{1,4}_{1,2}$ (see~\ref{eq:mulDisp}).
Accordingly we set:
\begin{equation}\label{eq:mysymgr}
\mathbb{G}^\textnormal{R}(\qq)=0, \qquad \textnormal{when $L$ is even, $d=2$, $q_1=q_2=\pi$},
\end{equation}
which enforces strain compatibility and stress equilibrium, in the sense of~(\ref{eq:discFull}).

\subsection{Three dimensions}
We follow the same methodology in 3D.
The equilibrium and strain admissibility conditions~(\ref{eq:discFour})
are unchanged, as well as the expression for the vectors $\kk^\textnormal{C,W}$
in~(\ref{eq:mul2}) and~(\ref{eq:wil072})
resulting from~(\ref{eq:mul1}) and~(\ref{eq:wil071}).
In 3D, we also extend~(\ref{eq:wil14}) as:
\begin{equation}\label{eq:wil14b}
k_i^\textnormal{R}(\qq)=\frac{\ii}{4}\tan\left(\frac{q_i}{2}\right)\left(1+\ee^{\ii q_1}\right)\left(1+\ee^{\ii q_2}\right)\left(1+\ee^{\ii q_3}\right),
\end{equation}
for the rotated scheme.
The strain and stress fields are now evaluated at the centers of the voxels
and the displacement field at their corners.
Derivatives of the displacement are estimated 
by differences at opposite corners.
For the strain components $\varepsilon_{11}$ and $\varepsilon_{12}$:
\begin{subequations}\label{eq:disc3deps}
\begin{eqnarray}\label{eq:disc3d11}
\varepsilon_{11}(\xx)&\approx& \frac{1}{4}
  \sum_{\substack{m=\pm 1\\ n=\pm 1}}\left[u_1\left(\xx+\frac{\eee_1+m\eee_2+n\eee_3}{2}\right)
                        -u_1\left(\xx-\frac{\eee_1-m\eee_2-n\eee_3}{2}\right)\right],\\
\varepsilon_{12}(\xx)&\approx& 
\frac{1}{8}
  \sum_{\substack{m=\pm 1\\ n=\pm 1}}\left[u_2\left(\xx+\frac{\eee_1+m\eee_2+n\eee_3}{2}\right)
                        -u_2\left(\xx-\frac{\eee_1-m\eee_2-n\eee_3}{2}\right) \nonumber\right.\\
  && +\left.
  u_1\left(\xx+\frac{\eee_2+m\eee_1+n\eee_3}{2}\right)
                        -u_1\left(\xx-\frac{\eee_2-m\eee_1-n\eee_3}{2}\right)\right], \label{eq:disc3d12}
\end{eqnarray}
\end{subequations}
where $\xx$ lie at the center of a voxel.
The expression for the strain component $\varepsilon_{22}$ (resp. $\varepsilon_{33}$) is obtained after
exchanging the indicia $1$ and $2$ (resp. $1$ and $3$) in~(\ref{eq:disc3d11}).
The component $\varepsilon_{23}$ (resp. $\varepsilon_{13}$) is derived
from~(\ref {eq:disc3d12}) by exchanging the indicia $3$ and $1$ (resp. $3$ and $2$).
Stress divergence is discretized in a similar manner.
Its first component reads:
\begin{eqnarray}\label{eq:disc3dsig}
\partial_i\sigma_{i1}(\xx)&\approx& 
  \sum_{\substack{m=\pm 1\\ n=\pm 1}}\left[\sigma_{11}\left(\xx+\frac{\eee_1+m\eee_2+n\eee_3}{2}\right)
                        -\sigma_{11}\left(\xx-\frac{\eee_1-m\eee_2-n\eee_3}{2}\right) \right.\nonumber\\
        && \left. +
         \sigma_{12}\left(\xx+\frac{\eee_2+m\eee_1+n\eee_3}{2}\right)
                        -\sigma_{12}\left(\xx-\frac{\eee_2-m\eee_1-n\eee_3}{2}\right) \right.\nonumber\\
        && \left. +
        \sigma_{13}\left(\xx+\frac{\eee_3+m\eee_1+n\eee_2}{2}\right)
                        -\sigma_{13}\left(\xx-\frac{\eee_3-m\eee_1-n\eee_2}{2}\right)\right],
\end{eqnarray}
where $\xx$ lie at one of the {\emph{edges}} of a voxel.
The components $\partial_i\sigma_{i2}$ and $\partial_i\sigma_{i3}$ 
are obtained from~(\ref{eq:disc3dsig}) by circular permutations of the indicia $1$, $2$ and $3$.
We note that~(\ref{eq:disc3deps}) and (\ref{eq:disc3dsig}) are the natural generalization of~(\ref{eq:discFull})
to $d=3$, expressed in the Cartesian basis $(\eee_1;\eee_2;\eee_3)$.

In 3D, Eq.~(\ref{eq:discFour2}) yields, 
for the Green operator:
\begin{equation}\label{eq:green3}
G'_{ij,kl}(\qq)=
\frac{
  \left(\lambda^0+2\mu^0\right) 
        \mysym{\left(r_ir_l^*\delta_{jk}\right)}
    +\lambda^0\left[    \mysym{\left(r_ir_l^* s_{jk}\right)}-\real\left(r_ir_j^*\right)\real\left(r_kr_l^*\right)\right]
   - \mu^0 r_ir_jr_k^*r_l^*
 }{
    \mu^0\left[2(\lambda^0+\mu^0)-\lambda^0\left|r_1^2+r_2^2+r_3^2\right|^2\right] 
 }
\end{equation}
where again $r_i=k_i/|\kk|$ and $\bm{s}$ is the symmetric second-order tensor:
\begin{equation}
s_{jj}=4\imag(r_ir_k^*)^2, \quad
s_{jk}=-4\imag(r_kr_j^*)\imag(r_kr_i^*), \quad 
i\neq j\neq k\neq i,
\end{equation}
with $\imag(\cdot)$ the imaginary part of the enclosed complex quantity.
Like in 2D, the operator $\mathbb{G}'$ follows minor and major symmetries~(\ref{eq:greenSym}).

Again, the operators  $\mathbb{G}^\textnormal{C}$,
$\mathbb{G}^\textnormal{W}$ and $\mathbb{G}^\textnormal{R}$
are derived using the expression for $\mathbb{G}'$
in~(\ref{eq:green3}) with $\kk=\kk^\textnormal{C}$, $\kk^\textnormal{W}$ and $\kk^\textnormal{R}$, respectively.
The symmetries~(\ref{eq:greenSym2}) and~(\ref{eq:grsym}) are verified in 3D as well.
But again, a special treatment is needed for $\mathbb{G}^\textnormal{C}$ and $\mathbb{G}^\textnormal{R}$
when $L$ is even, at the modes $\qq$ for which $\kk^\textnormal{C,R}(\qq)=0$.
Like in 2D, the displacement is undefined at these frequencies and the Fourier coefficients
of the strain field are zero and so
we set $\mathbb{G}^\textnormal{C,R}=0$ at these frequencies.
More precisely, we apply (\ref{eq:mulfix}) when $d=3$ and,
for the rotated scheme:
\begin{equation}
\mathbb{G}^\textnormal{R}(\qq)=0, \qquad \textnormal{ if $L$ is even, $d=3$ and $q_i=q_j=\pi$ with $i\neq j$}.
\end{equation}

The operators $\mathbb{G}^\textnormal{C}$, $\mathbb{G}^\textnormal{W}$ and $\mathbb{G}^\textnormal{R}$ are,
in 2D and 3D, periodic functions where, 
contrary to $\mathbb{G}$, high frequencies are cut.
Accordingly, we expect faster convergence rates for schemes using operators derived from finite differences
and more exact local fields,
     as was previously observed in the conductivity problem~\cite{Willot14a}.
     We also expect higher accuracy for the local fields when employing $\mathbb{G}^\textnormal{R}$
     rather than the other discrete operators $\mathbb{G}^\textnormal{W}$ and $\mathbb{G}^\textnormal{C}$.
     First, the operator $\mathbb{G}^\textnormal{R}$ is based on centered differences
     which are more precise than forward and backward differences, used in $\mathbb{G}^\textnormal{W}$.
     Second, derivatives are evaluated more locally when using $\mathbb{G}^\textnormal{R}$
     rather than $\mathbb{G}^\textnormal{C}$.
     Indeed, the latter are computed at points 
     separated by $2$ voxels for $\mathbb{G}^\textnormal{C}$ instead of $\sqrt{2}$ (in 2D) 
     or $\sqrt{3}$ voxels (in 3D) for $\mathbb{G}^\textnormal{R}$.
     
The above considerations guided
the choice for the discretization
schemes~(\ref{eq:disc3deps})
and~(\ref{eq:disc3dsig}), leading to
$\kk^\textnormal{R}$ and $\mathbb{G}^\textnormal{R}$.
Clearly, many other choices are possible,
and Eq.~(\ref{eq:green3})
gives a general class of Green operators based on finite-differences.
The latter depend on the choice for the complex vector $\kk$.
However, a systematic investigation of such discrete schemes
is beyond the present study.

In the the rest of this study, we estimate the accuracy
of the local fields and of the effective properties 
predicted by the various schemes, as well as
their convergence rates.
We denote by DS and AS the direct and accelerated schemes 
defined by~(\ref{eq:ds}) and~(\ref{eq:as}) respectively,
when $\mathbb{G}$ is used.
We denote by
$DS_\textnormal{C}$, $DS_\textnormal{W}$, $DS_\textnormal{R}$
and $AS_\textnormal{C}$, $AS_\textnormal{W}$ and $AS_\textnormal{R}$,
the same algorithms obtained by substituting $\mathbb{G}$ with $\mathbb{G}^\textnormal{C,W,R}$ respectively.
We emphasize that, for a given
Green operator, the direct and accelerated schemes produce the same strain and stress fields,
up to round-off errors.

\section{Local strain and stress fields accuracy}\label{sec:accuracy}
\subsection{Two-dimensions}\label{sec:2dacc}
\begin{figure}
\begin{center}
\includegraphics[width=4cm]{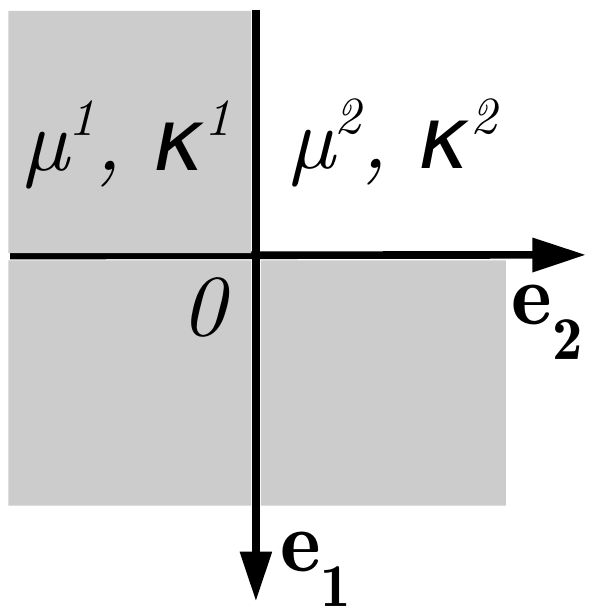}
\caption{\label{fig:micro2d}
Elementary periodic domain $\Omega=[−1/2;1/2]^2$
containing a square inclusion with elastic moduli $\mu^2$, $\kappa^2$ (top-left, shown in white) embedded in a matrix
(shown in gray) with elastic moduli $\mu^1$, $\kappa^1$ .
}
\end{center}
\end{figure}

Hereafter we consider the 2D `four-cell' microstructure,
where the periodic domain $\Omega$ is divided into $4$ identical squares
of surface fraction $25$\%. 
Its non-trivial solution with singular fields at the corners
makes it a good benchmark for numerical schemes.
Furthermore, the microstructure
is discretized exactly at any resolution, provided $L$ is even.
In the following, we make use of a simplified version of the four-cell microstructure
made of a single quasi-rigid square inclusion embedded in a matrix
(Fig.~\ref{fig:micro2d}). We set the contrast to $\chi=10^3$.
The material is subjected to the macroscopic strain loading:
\begin{equation}\label{eq:load2d}
  \overline{\varepsilon_{ij}}=\frac{1}{2}\left(\delta_{i1}\delta_{j2}+\delta_{j1}\delta_{i2}\right).
\end{equation}
We determine the strain and stress fields predicted by FFT schemes when using the Green
operators $\mathbb{G}$ or $\mathbb{G}^\textnormal{C,W,R}$.
The fields are computed using the accelerated scheme~(\ref{eq:as})
at discretizations $L=512$, $1024$ and $2048$.
Iterations are stopped when the strain and stress fields maximum variation over two iterations in any pixel is 
less than $2\,10^{-13}$. These variations are the effect of round-off errors in double precision floating point numbers.
These computations allow us to compare the effect of the discretization, independently of the algorithm used for convergence.

\begin{figure}
\begin{center}
\begin{tabular}{p{0.3cm}cccc}
  & $L=512$ & $L=1024$ & $L=2048$ \\
\vspace{\lenx}$\mathbb{G}$ &
\includegraphics[width=\len]{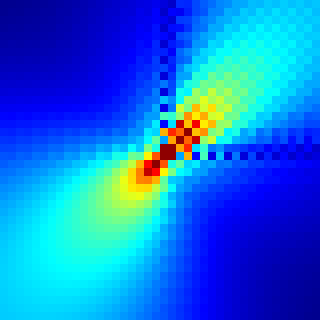} &
\includegraphics[width=\len]{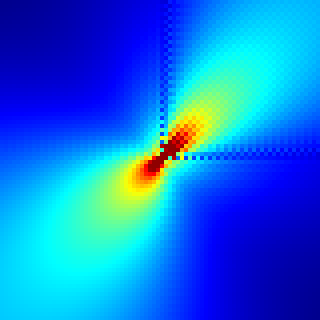} &
\includegraphics[width=\len]{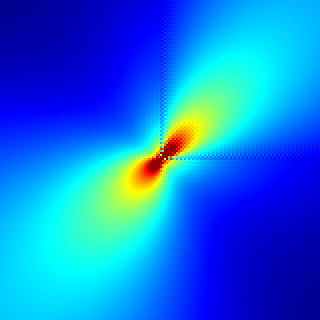} \\
\vspace{\lenx}$\mathbb{G}^\textnormal{C}$ &
\includegraphics[width=\len]{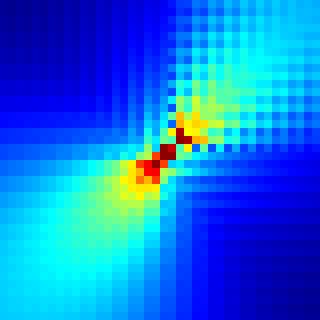}  &
\includegraphics[width=\len]{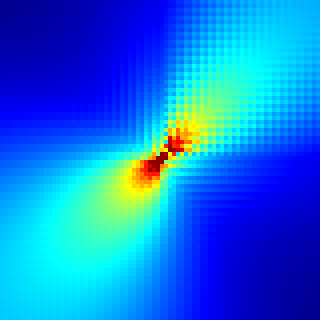}  &
\includegraphics[width=\len]{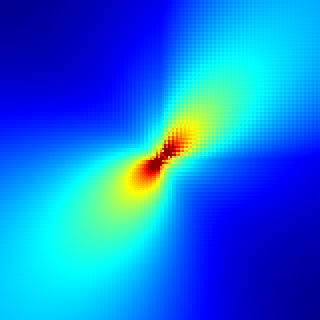}  \\
\vspace{\lenx}$\mathbb{G}^\textnormal{W}$ &
\includegraphics[width=\len]{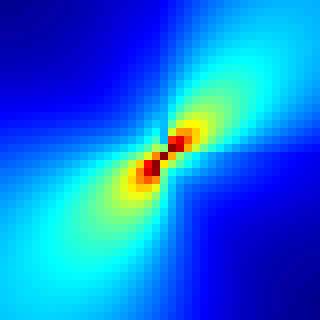}    &
\includegraphics[width=\len]{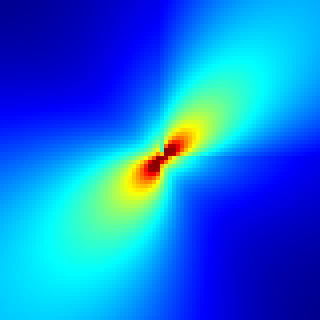}    &
\includegraphics[width=\len]{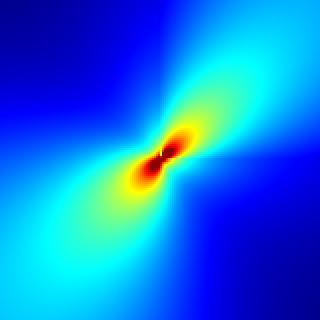}    \\
\vspace{\lenx}$\mathbb{G}^\textnormal{R}$ &
\includegraphics[width=\len]{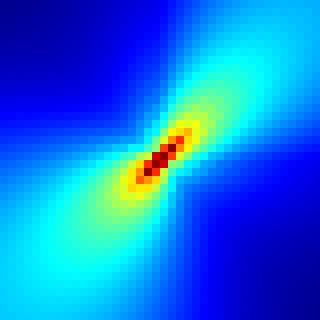}    &
\includegraphics[width=\len]{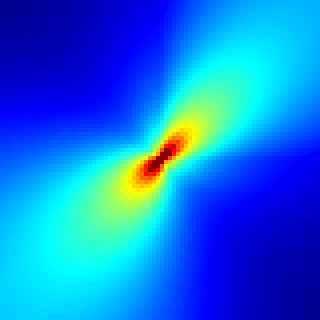}    &
\includegraphics[width=\len]{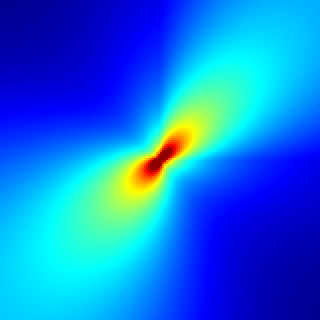}    
\end{tabular}
\caption{\label{fig:accuracy2d}
Stress component $\sigma_{12}(\xx)$ 
predicted by the various FFT schemes at the three resolutions $L=512$, $1024$ and $2048$ (left to right)
in the region $[-0.04;0.04]^2$.
The center of the region is the bottom-left corner 
of the square inclusion in Fig.~(\ref{fig:micro2d}).
}
\end{center}
\end{figure}

We focus on the stress component $\sigma_{12}(\xx)$ parallel to the applied
loading in a small region $[-0.04; 0.04]^2$ around the corner of the inclusion (Fig.~\ref{fig:accuracy2d}).
At low resolution $L=512$, numerical methods
predict values as large as $10.1$ in a few pixels, because of the singularity of the stress field at the corner.
To highlight the field patterns,
we threshold out the values above $3.5$, which amount to $0.24$\% of the pixels.
Using the same color scale for all images, the smallest stress value, equal to $1.5$, is
shown in dark blue whereas the highest, equal to $3.5$, is in dark red. Green, yellow and orange lie in-between.

As expected, in the limit of very fine resolution, all methods tend to the same local stress field,
as shown by the similar field maps obtained at resolution $L=2048$. 
However, use of the Green operator $\mathbb{G}$ leads to
spurious oscillations along the interfaces of the inclusion, 
up to resolutions as big as $2048^2$ pixels, a side-effect noticed in~\cite{Willot14a} in conductivity.
The oscillations do not disappear after computing local averages of the fields (not shown).

Strong oscillations are produced by schemes using~$\mathbb{G}^\textnormal{C}$ as well, not only
in the quasi-rigid inclusion, but also in the matrix.
We observe checkerboard patterns in the former, and vertical and horizontal alignments in the latter,
at resolution $1024^2$.
These oscillations are greatly reduced by the use of $\mathbb{G}^\textnormal{W}$.
Still, due to the non-symmetric nature of $\mathbb{G}^\textnormal{W}$,
the stress is not correctly estimated along a 
line of width $1$ pixel oriented upward from the inclusion corner.
Similar patterns are observed, in other directions, along the three other corners of the inclusion (not shown).
These issues are solved when using $\mathbb{G}^\textnormal{R}$ which produces a stress field that respects the symmetries of the problem.
Furthermore, use of $\mathbb{G}^\textnormal{R}$ greatly reduces oscillations compared to $\mathbb{G}$
and $\mathbb{G}^\textnormal{C}$.

\subsection{Three-dimensions}
In this section, we consider a 3D material analogous to the four-cell microstructure in 2D.
We divide the periodic domain into $8$ identical cubes of volume fraction
$12.5$\%. One is the inclusion,
the other $7$ are the matrix.
To highlight the symmetries of the problem, we assume the inclusion is centered in the domain $\Omega$
and contained in the region $[-1/4;1/4]^3$.
Again, we apply a macroscopic strain loading of the form~(\ref{eq:load2d}).
The inclusion is quasi-rigid compared to the matrix with contrast of properties $\chi=10^3$.
We compute the strain and stress fields predicted by each Green operator using the accelerated scheme.
As in Sec.~(\ref{sec:2dacc}) we let the iterations converge up to round-off
errors in double precision.

\begin{figure}
\begin{center}
\begin{tabular}{p{0.3cm}ccc}
  & $L=256$ & $L=512$ & $L=1024$ \\
\vspace{\lemx}$\mathbb{G}$ &
\includegraphics[width=\lem]{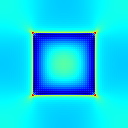} &
\includegraphics[width=\lem]{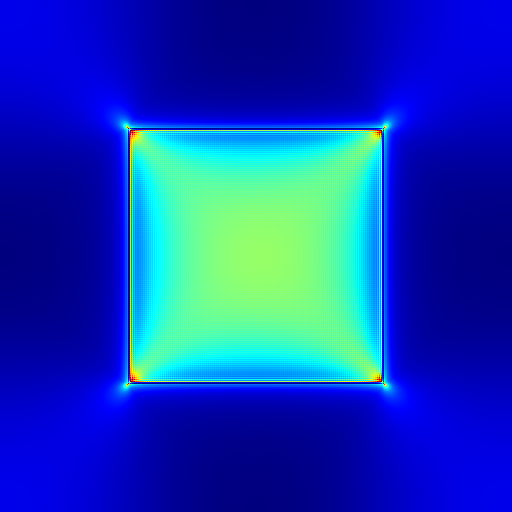} &
\includegraphics[width=\lem]{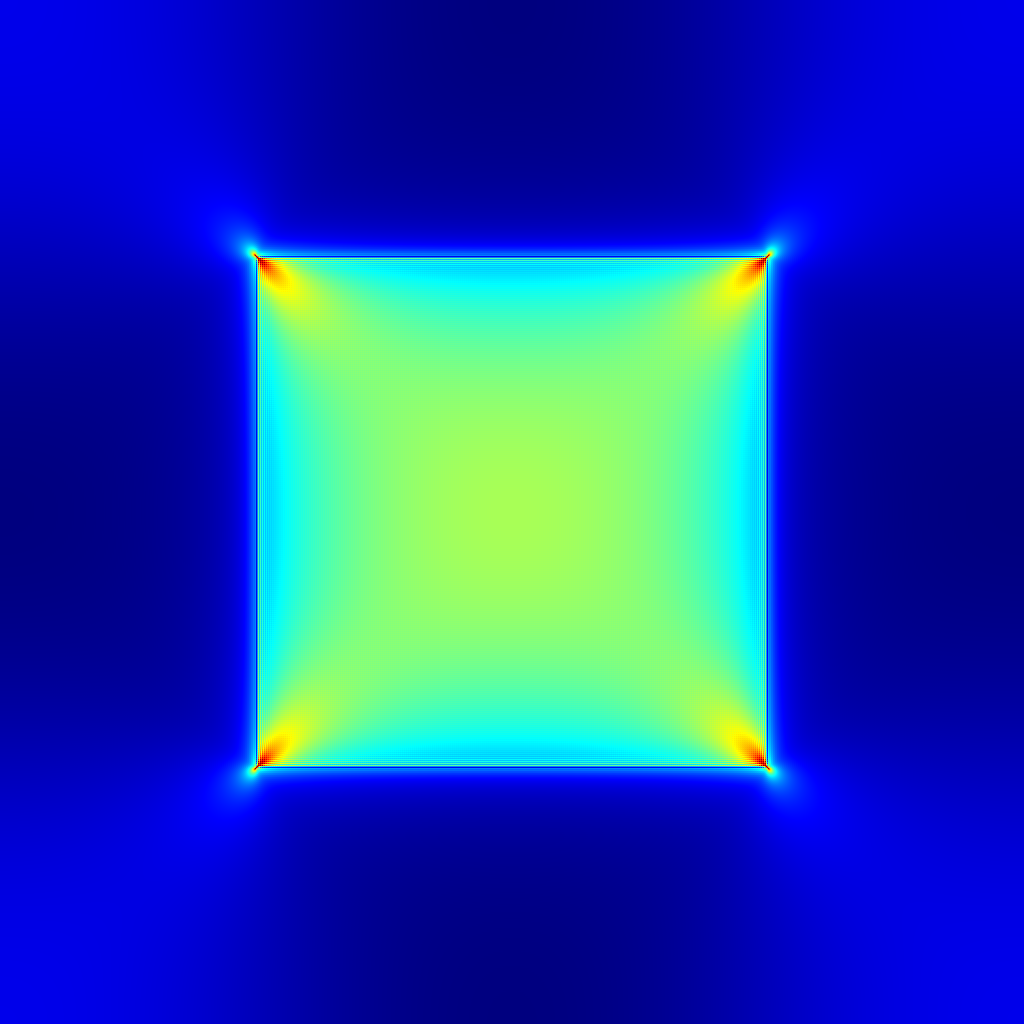} \\
\vspace{\lemx}$\mathbb{G}^\textnormal{C}$ &
\includegraphics[width=\lem]{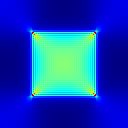}  &
\includegraphics[width=\lem]{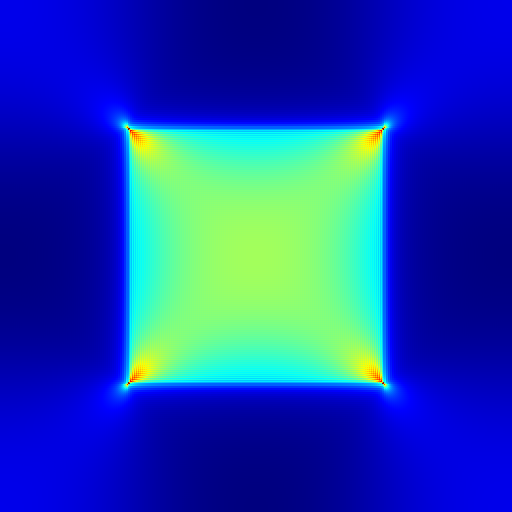}  &
\includegraphics[width=\lem]{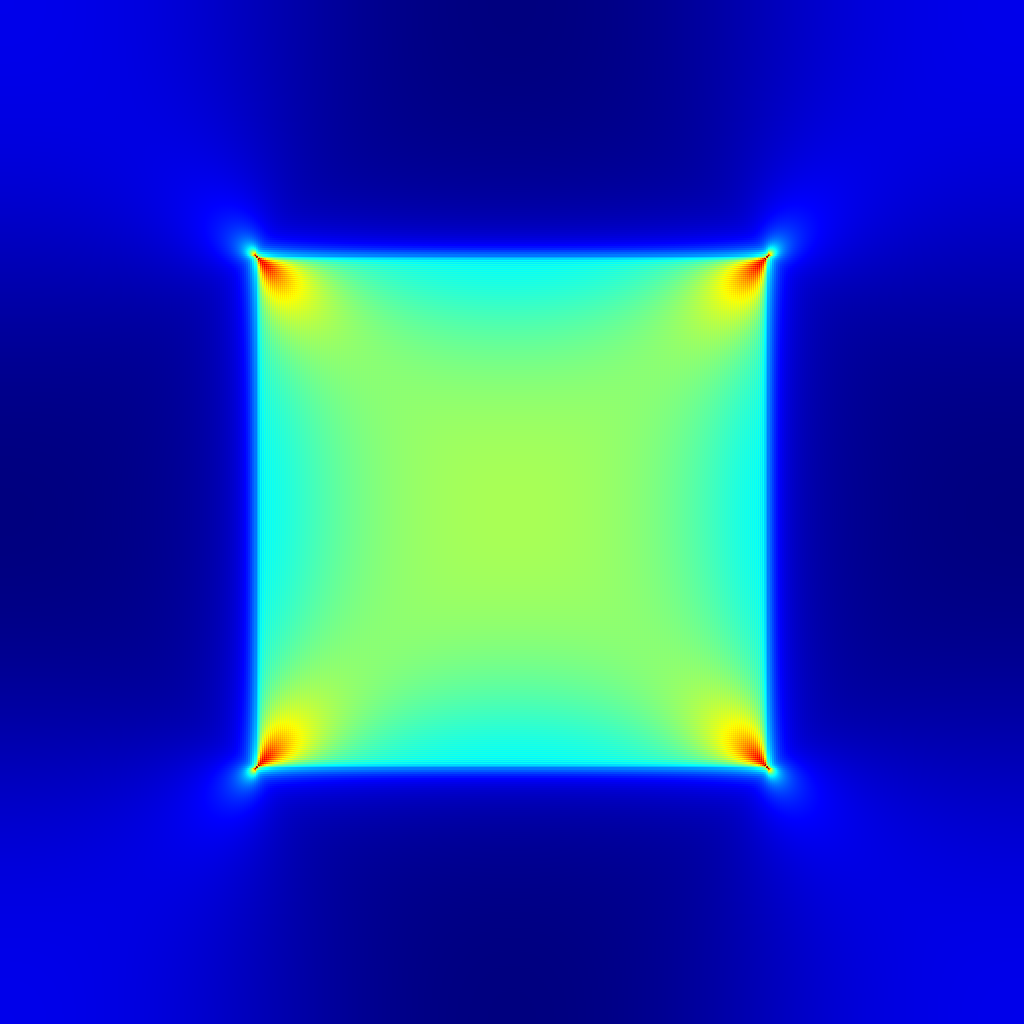}  \\
\vspace{\lemx}$\mathbb{G}^\textnormal{W}$ &
\includegraphics[width=\lem]{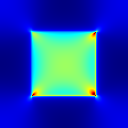}    &
\includegraphics[width=\lem]{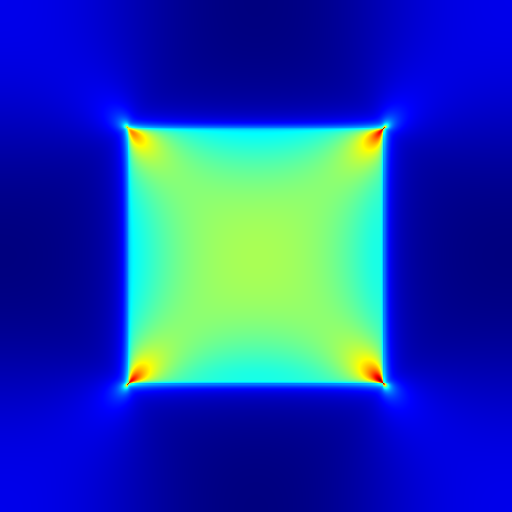}   & 
\includegraphics[width=\lem]{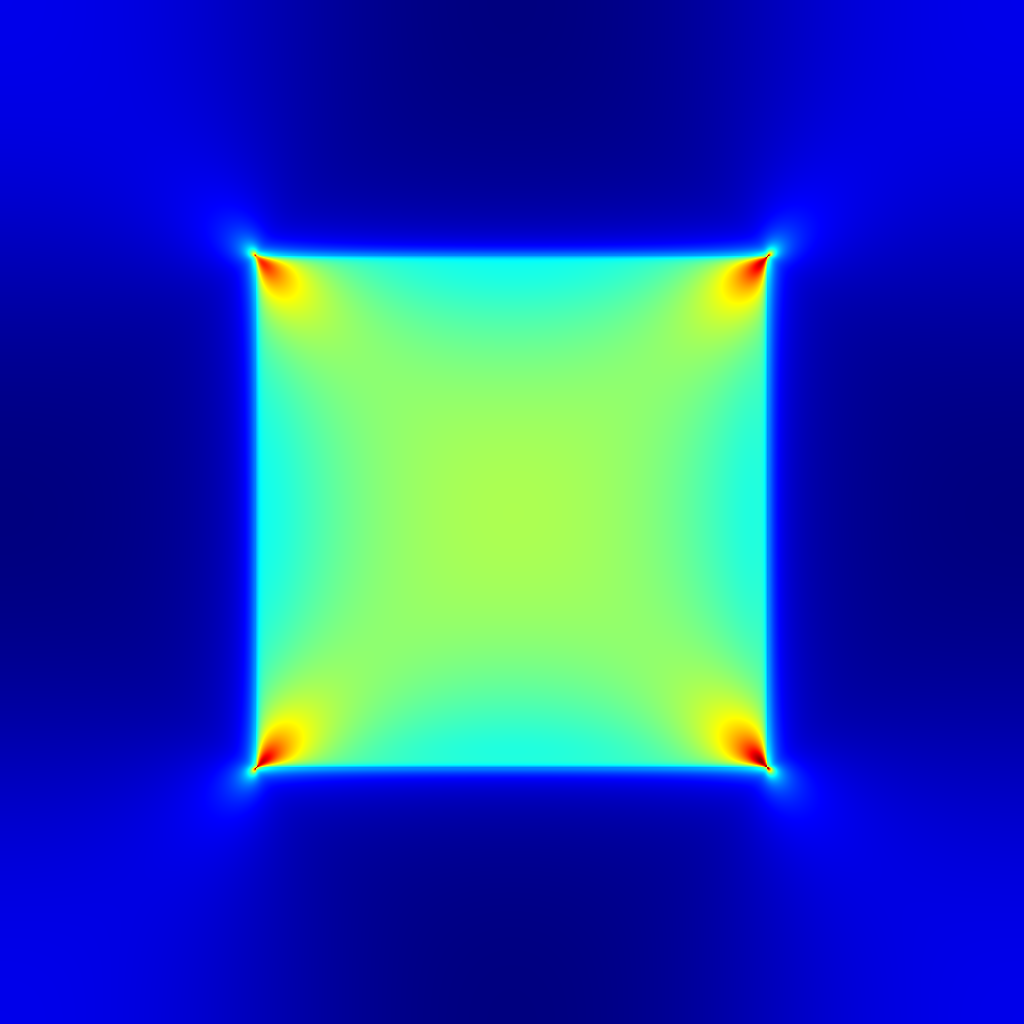}   \\
\vspace{\lemx}$\mathbb{G}^\textnormal{R}$ &
\includegraphics[width=\lem]{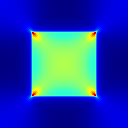}    &
\includegraphics[width=\lem]{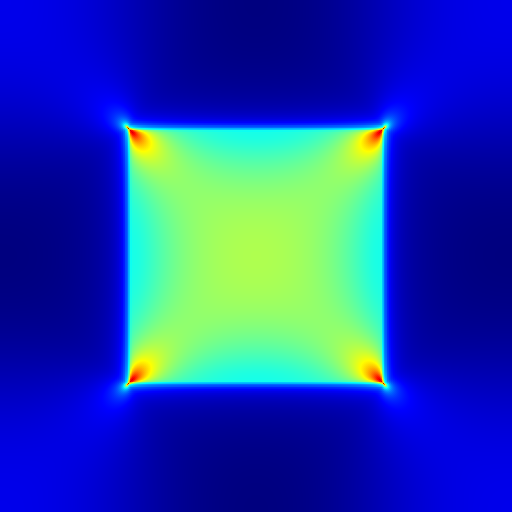}   &  
\includegraphics[width=\lem]{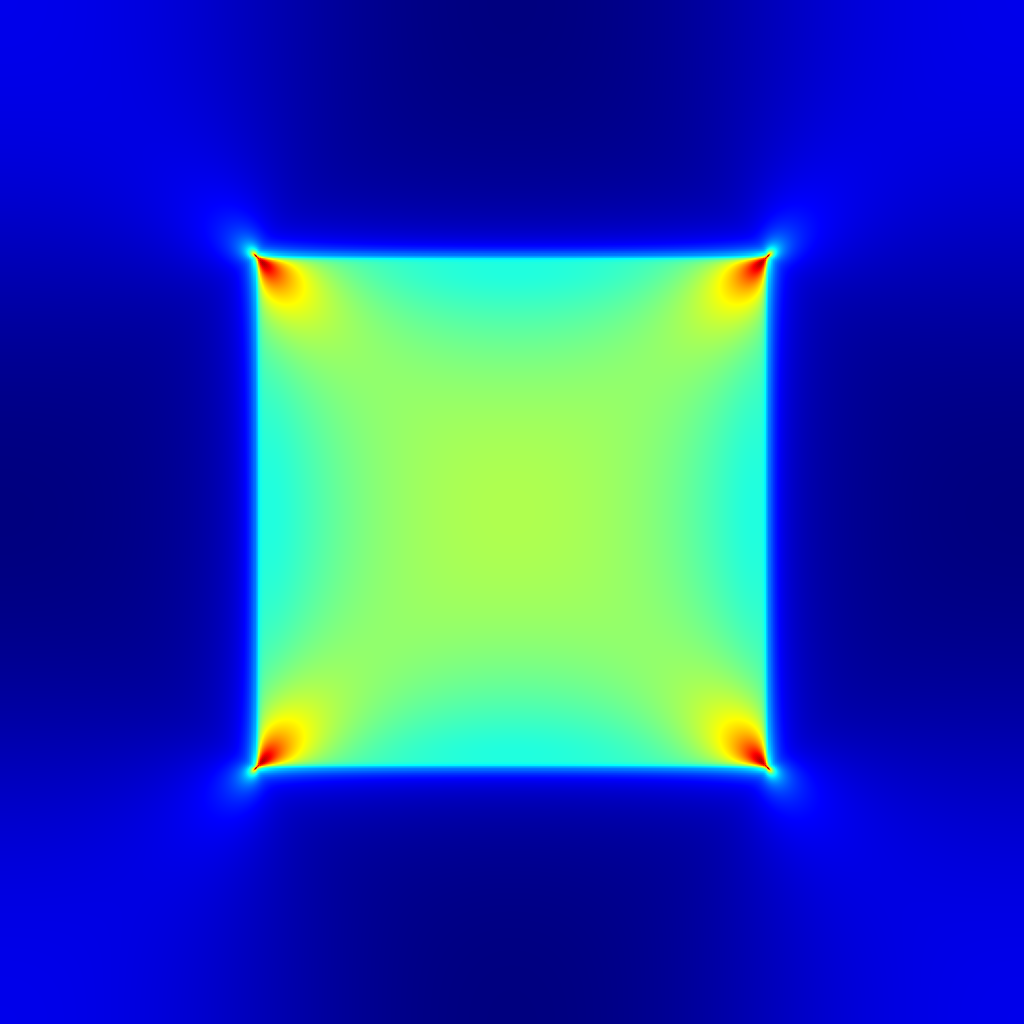}    \\
\end{tabular}
\caption{\label{fig:accuracy3d}
2D section of the stress component $\sigma_{12}(\xx)$ along the plane $x_3=-0.2461$
predicted by the various FFT schemes at the three resolutions $L=256$, $L=512$ and $L=1024$ (left to right).
The section is parallel to one of the faces of the inclusion and close to the interface with the matrix.
}
\end{center}
\end{figure}

A 2D section of the stress component $\sigma_{12}(\xx)$ is represented in Fig.~(\ref{fig:accuracy3d}).
The section is a cut parallel to one of the faces of the inclusion, normal to $\eee_3$,
of equation $x_3=-0.2461$. The section intersects the inclusion, but is very close to the interface with the matrix.
Again, to highlight the field
patterns, we threshold out values of the field greater than $8.5$, this time less than $0.04$\% of the voxels,
and represent all field using the same color scale.

At high resolution $L=1024$, the fields 
resulting from the use of $\mathbb{G}$ and $\mathbb{G}^\textnormal{C,W,R}$ are close to one another.
However, stress patterns near the corners of the inclusion are less pronounced with $\mathbb{G}$ than with the other methods.
At smaller resolutions $L=256$ and $L=512$, the stress fields predicted by $\mathbb{G}$
are notably different from the others,
suggesting slower size-convergence with this operator.
Furthermore, the field maps computed at resolution $L=512$
confirms the results obtained in 2D:
strong oscillations are observed inside the inclusion
when using $\mathbb{G}$ and $\mathbb{G}^\textnormal{C}$.
The two methods produce artificial 
patterns directed vertically and horizontally,
close to the interface.
Conversely, the fields produced by $\mathbb{G}^\textnormal{W}$ and $\mathbb{G}^\textnormal{R}$
have the smallest oscillations,
but that of $\mathbb{G}^\textnormal{W}$ are not symmetric.
When $L=256$, indeed, the stress field near the top-left corner
of the inclusion stands out from that in the other corners.
This effect only slowly disappears when $L$ is increased.
The solution resulting from the use of $\mathbb{G}^\textnormal{R}$ does not suffer from this  problem.
As in 2D, it
produces symmetric fields. Furthermore, the latter are close to one another at all resolutions and contain 
almost no oscillations.

\subsection{Periodic array of spheres}
\begin{figure}
\begin{center}
\includegraphics[width=9cm]{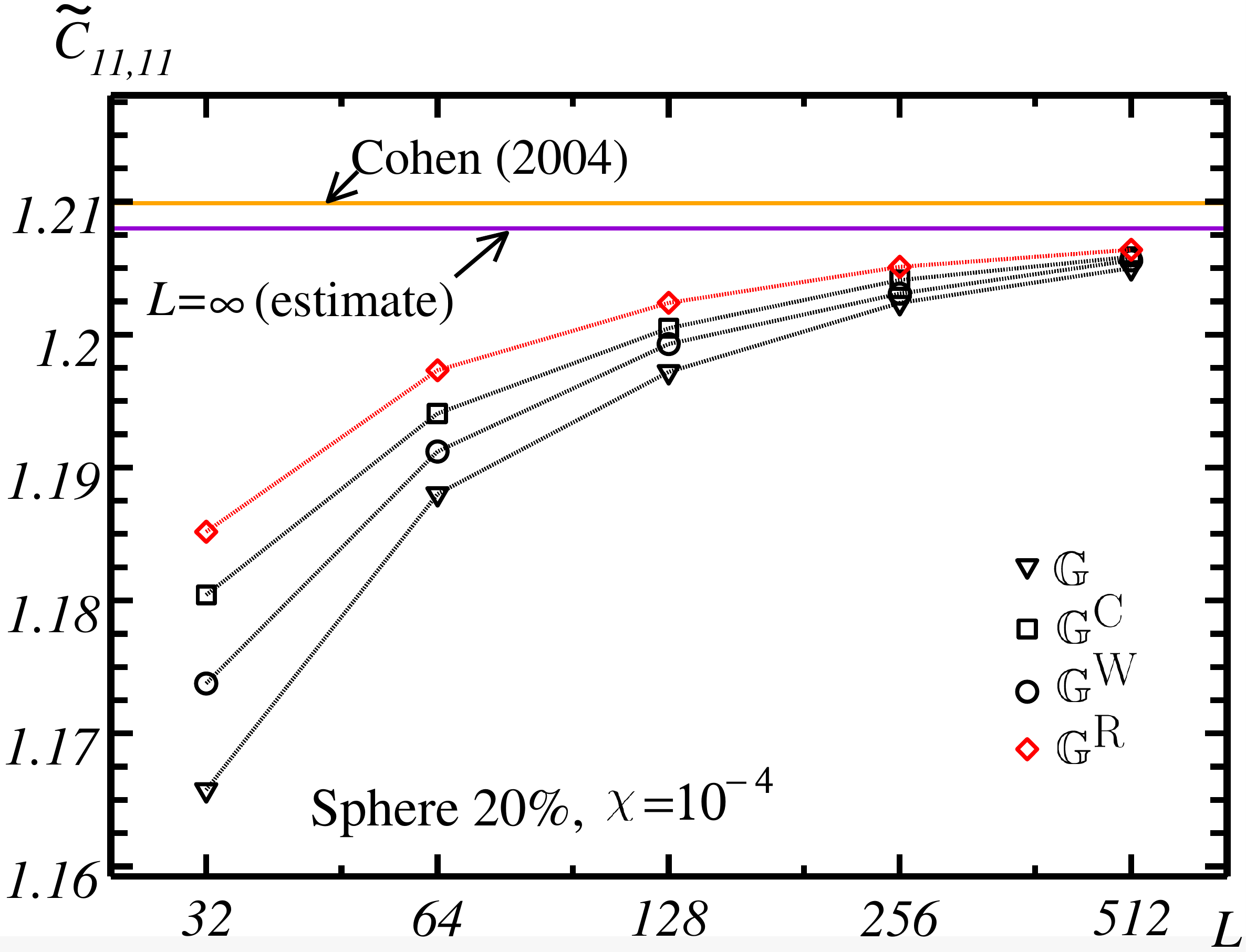}
\caption{\label{fig:eff}
Apparent elastic modulus $\widetilde{C}_{11,11}$ estimated 
by FFT methods using the Green operators $\mathbb{G}$ and
$\mathbb{G}^\textnormal{C,W,R}$ (black and red),
at increasing
resolution $L$.
Orange: estimate in~\cite{Cohen04};
violet: estimate of the asymptotic effective modulus $\widetilde{C}_{11,11}$ using FFT data.
}
\end{center}
\end{figure}

Contrary to the previous sections, we now consider a microstructure without singularities (edges or corners)
and focus on the effect of the Green operator discretization on the effective elastic properties.
In the rest of this section, the elementary domain $\Omega$ contains one spherical inclusion
of volume fraction $20$\%, so that the material is a
periodic array of spheres. The spheres are very soft with contrast of properties $\chi=10^{-4}$.
We compute the effective elastic modulus $\widetilde{C}_{11,11}$
produced by either $\mathbb{G}$ or
$\mathbb{G}^\textnormal{C,W,R}$ at increasing resolutions $L=32$, $64$, $128$, $256$ and $512$.
Again, we use the accelerated scheme and
iterations are stopped when the stress field
maximum variation over two iterations in any pixel is 
less than $2\,10^{-10}$.
Results are shown in Fig.~(\ref{fig:eff})
and are compared with the analytical estimate in~\cite{Cohen04}.
When the resolution increases, 
the effective elastic modulus $\widetilde{C}_{11,11}$
increases up to a limit value that we estimate to about $1.208\pm 0.001$, for all schemes.
As observed in other studies~\cite{Dunant13},
very large systems are needed to compute this estimate at a high precision.

This is especially true 
of the Green operator $\mathbb{G}$ which has the slowest convergence with respect to the system size.
At fixed resolution, the error on the predictions given by $\mathbb{G}$
is about $2$ times larger than the one provided by $\mathbb{G}^\textnormal{R}$,
which, among all methods, gives the best estimate.
The operators $\mathbb{G}^\textnormal{C}$ and $\mathbb{G}^\textnormal{W}$ stand in-between.
This is another indication of the benefits of 
the operator $\mathbb{G}^\textnormal{R}$.

\section{Convergence rate}\label{sec:rate}
\subsection{Convergence rate with respect to stress equilibrium}
In this section, we estimate the rates of convergence of the direct and accelerated schemes
DS, DS$_\textnormal{C,W,R}$, AS and AS$_\textnormal{C,W,R}$,
that use the various Green operators.
All schemes enforce stress equilibrium at convergence only,
therefore we follow~\cite{moulinec13} and consider a criterion
based on the  $L^2$-norm:
\begin{equation}
\label{eq:criterion}
\eta=\frac{||\textnormal{div}(\sig)||_2}{||\langle\sig\rangle||}=\frac{1}{||\langle\sig\rangle||}
\sqrt{\frac{1}{|\Omega|}\int_\Omega{\rm d}^d\xx\,|\textnormal{div}(\sig)|^2}=
\frac{1}{||\langle\sig\rangle||}
\sqrt{\sum_{\qq}\left|\kk(\qq)\cdot\sig(\qq)\right|^2},
\end{equation}
where $\eta\ll 1$ is the precision and the normalizing factor $||\langle\sig\rangle||$ is the Frobenius norm:
$$
||\langle\sig\rangle||^2=\sum_{i, j}\langle\sigma_{ij}(\xx)\rangle^2.
$$
In~(\ref{eq:criterion}) we set $\kk=\kk^\textnormal{C,W,R}$ for the schemes
using $\mathbb{G}^\textnormal{C,W,R}$ and $\kk=\ii\qq$ when using the Green operator $\mathbb{G}$,
so that $\kk\cdot\sig(\qq)$ is
the divergence of the stress field in the Fourier domain, estimated according to the various discretization schemes.

We now estimate the convergence rates on a random microstructure.
In the following, the domain $\Omega$ is a
(periodized) Boolean model of spheres of resolution $L=64$ and volume fraction $17$\%,
below the percolation threshold of the spheres --- of about $29$\%~\cite{willot09}.
To obtain meaningful comparisons, we use the same randomly-generated microstructure
for all schemes. This particular 
configuration contains $743$ spheres of diameter $5$ voxels, about $13$ times smaller than the size of $\Omega$.

\begin{figure}
\begin{center}
\includegraphics[width=9cm]{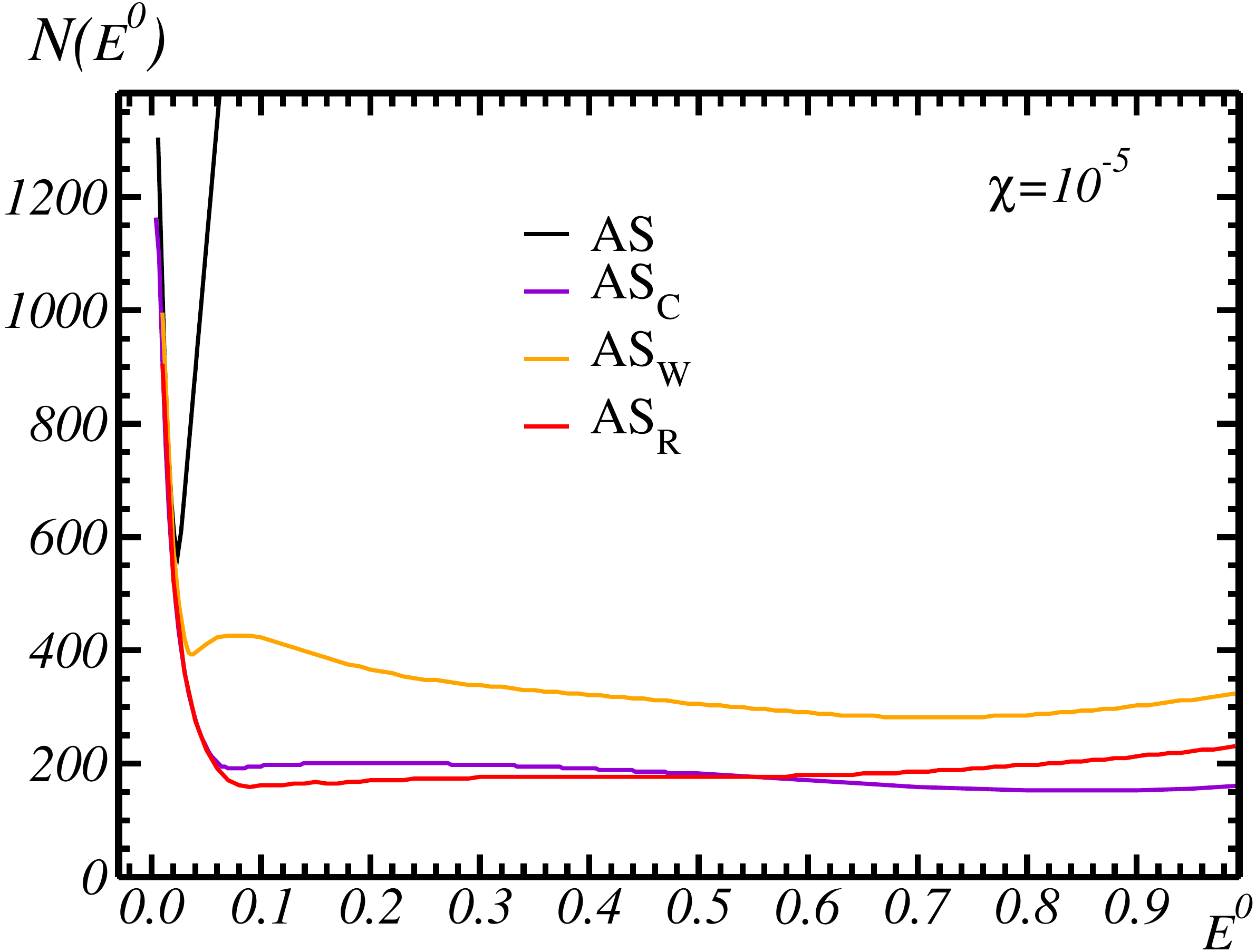} 
\caption{\label{fig:convref0}
Number of iterations $N(E^0)$ 
required to achieve convergence,
as a function of the reference Young modulus $E^0$,
for the accelerated schemes AS and AS$_\textnormal{C,W,R}$ using various Green operators.
Convergence is achieved when the precision $\eta=10^{-8}$ is reached.
The microsctructure is a Boolean model of quasi-porous spheres with $\chi=10^{-5}$.
}
\end{center}
\end{figure}

Taking $\nu^0=\nu^1=\nu^2=0.25$ for the reference Poisson ratio, we
compute numerically the number of iterations $N(E^0)$
required to reach the precision $\eta\leq 10^{-8}$, for varying
reference Young moduli $E^0$, in the range $0<E^0<1$.
We consider the Boolean model of spheres with contrast $\chi=10^{-5}$ 
and the various accelerated schemes AS and AS$_\textnormal{C,W,R}$ (Fig.~\ref{fig:convref0}).
Within the range $0<E^0\lesssim 0.03$,
the number of iterations $N(E^0)$ is about the same 
for all accelerated schemes.
When $E^0> 0.03$, however, 
$N(E^0)$ is a strongly increasing function of $E^0$ for scheme AS,
contrary to the other schemes AS$_\textnormal{C,W,R}$.
For the latter, $N(E^ 0)$ decreases with $E^0$ up to a local minimum,
beyond which variations are much less sensitive to $E^0$.
One unique local minimum around $E^0\approx 0.09$ is found for 
scheme  AS$_\textnormal{R}$, whereas the schemes AS$_\textnormal{W}$ and AS$_\textnormal{C}$
exhibit two local minima.

The effect of the Poisson ratio
is also investigated numerically.
We let $\nu^0=0.25\pm 0.01$ and $0.25\pm 0.05$ for various values of $E^0$ 
with $\chi=10^{-5}$
and observe a strong increase of the number of iterations $N(E^0)$,
for the schemes AS and AS$_\textnormal{C,W,R}$.
The same behavior is observed for the direct scheme DS and DS$_\textnormal{C,W,R}$
with $\chi=10^{-2}$.
Therefore, in the following, we fix the Poisson ratio to $\nu^0=0.25$ for the reference tensor, for all schemes
and all contrast of properties $\chi$.
This leaves one parameter, $E^0$, to optimize on.
We use the gradient descent method to determine
 a local minimum of $N(E^0)$
for arbitrary contrast and scheme DS, DS$_\textnormal{C,W,R}$, AS and AS$_\textnormal{C,W,R}$.
As above, $N(E^0)$ is the number of iterations necessary to reach $\eta\leq 10^{-8}$.
We choose $E^0=0.51(E^1+E^2)$ for schemes DS, DS$_\textnormal{C,W,R}$
and $E^0=\sqrt{E^1E^2}$ for schemes AS and AS$_\textnormal{C,W,R}$
as initial guess for $E^0$, suggested by~(\ref{eq:ref}) and~(\ref{eq:refAcc}).
At each step, we determine if $E^0$ is to be increased or decreased,
by comparing $N(E^0)$ with $N(E^0+\delta E^0)$
where $\delta E^0=0.01 E^0$. 
It frequently happens that $N(E^0)=N(E^0+\delta E^0)$.
In that case, we compare the values of the precision $\eta$ after $N(E^0)$ iterations
and follow the direction that minimizes $\eta$.
Iterations are stopped whenever $N(E^0)$ is unchanged
after two descent steps.

The gradient descent method determines a local minimum rather than the global minimum,
which is sub-optimal.
To check the validity of the results, further numerical investigations are carried out
for $\chi=10^{-2}$, $10^2$ and schemes DS and DS$_\textnormal{C,W,R}$.
The method predicts the global minimum in these cases.
This also holds for schemes AS and AS$_\textnormal{W,R}$
with $\chi=10^{-5}$, $10^5$,
but not for scheme $AS_\textnormal{C}$ with $\chi=10^{-5}$.
However, in this case the number of iterations $N(E^0)$ are very similar at the two 
local minima, as shown in Fig.~(\ref{fig:convref0}).
In the following, the results given by the gradient descent method 
are used as-is.

\begin{figure}
\begin{center}
\includegraphics[width=9cm]{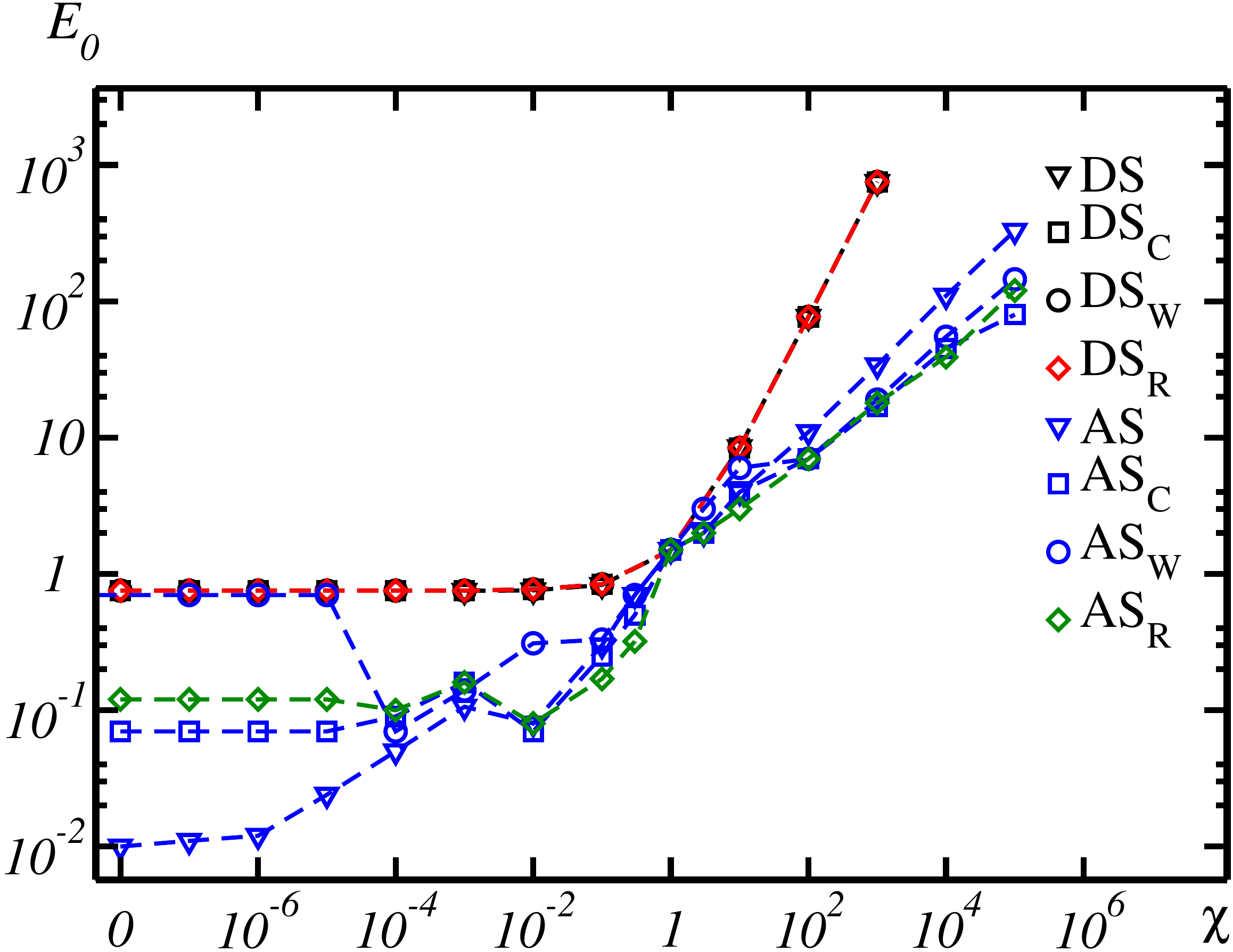} 
\caption{\label{fig:convref}
Optimal reference Young modulus $E^0$ as a function of the contrast of properties $\chi$,
for the various FFT methods, in log-log scale.
Direct schemes: black and red (nearly superimposed to one another);
accelerated schemes: blue and green.
Results for a porous material ($\chi=0$) are indicated at the left of the graph.
The material is a Boolean model of spheres with volume fraction $17$\%.
}
\end{center}
\end{figure}

Results for the optimal reference $E^0$ are indicated in Fig.~(\ref{fig:convref}).
For the direct scheme, the optimal reference follows~(\ref{eq:ref}) with 
$0.5003\leq\beta\leq 0.509$, independently of the Green operator used. Values of $\beta$ smaller than $1/2$ 
lead to non-converging schemes.
For the accelerated schemes, the situation is less simple, and differs depending on the Green operator in use.
For the scheme AS with Green operator~$\mathbb{G}$,
the choice~(\ref{eq:refAcc}) is optimal except in the region $\chi\leq 10^{-3}$
where the value of $E^0$ tend to a small constant of about $0.01$.
Similar behavior is found for the schemes AS$_\textnormal{C,R}$ for which:
\begin{equation}
   E^0\approx E^0_1+\sqrt{\chi},
\end{equation}
with $E^0_1=0.07$ for $AS_\textnormal{C}$
and $E^0_1=0.12$ for $AS_\textnormal{R}$.
Similar behavior has been observed numerically in~\cite{Willot14a}, in the context of conductivity.
For the scheme $AS_\textnormal{W}$ that uses $\mathbb{G}^\textnormal{W}$, 
the optimal choice for $E^0$ follows the same pattern as above with $E^0_1=0.7$
except in the region $10^{-4}\leq \chi\leq 10^{-1}$.
This behavior is an effect of 
the presence of two local minima, similar to that shown in Fig.~(\ref{fig:convref0}) for $\chi=10^{-5}$.

\begin{figure}
\begin{center}
\includegraphics[width=9cm]{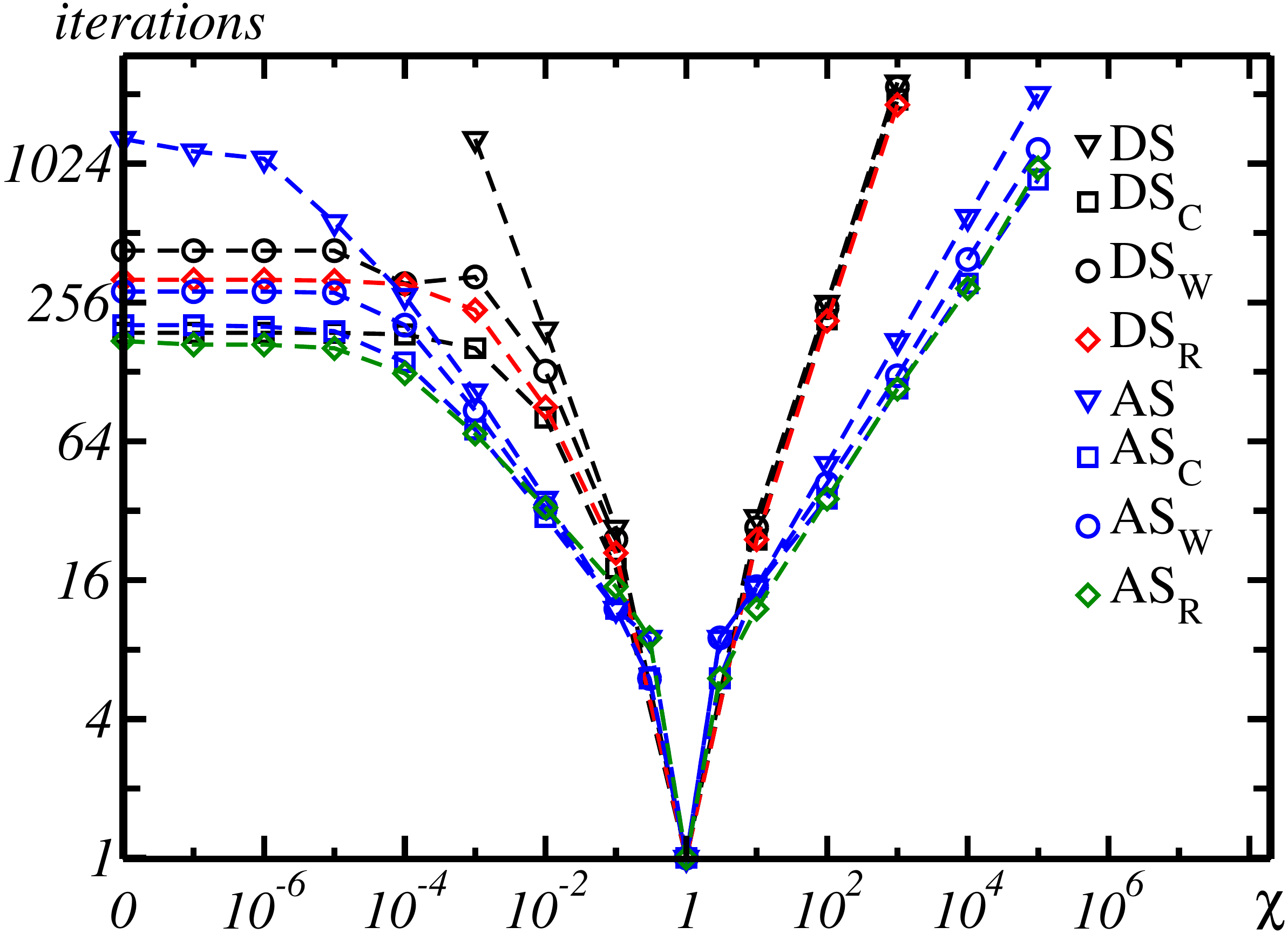} 
\caption{\label{fig:convrate}
Number of iterations as a function of the contrast of properties $\chi$,
for the various FFT methods, in log-log scale.
Direct schemes: black and red;
accelerated schemes: blue and green.
Results for a porous material ($\chi=0$) are indicated at the left of the graph.
The material is a Boolean model of spheres with volume fraction $17$\%.
}
\end{center}
\end{figure}

Convergence rates, computed with optimized reference, are represented in Fig.~(\ref{fig:convrate})
as a function of the contrast, in log-log scale. Results for $\chi=0$ (strictly porous media)
have been included in the same graph (left point).
As is well-established~\cite{Michel01,Moulinec03},
the number of iterations in the direct scheme DS scales as $\chi$ when
$\chi\gg 1$ and $1/\chi$ when $\chi\ll 1$.
For the accelerated scheme AS, the number of iterations is smaller and
follows $\sqrt{\chi}$ when
$\chi\gg 1$ and $1/\sqrt{\chi}$ when $\chi\ll 1$,
with one exception.
At very high contrast of properties $\chi<10^{-6}$, including at $\chi=0$,
convergence is reached after a finite number of iterations, about $1,300$.
This particular behavior is presumably sensitive to the value choosen for the
requested precision $\eta=10^{-8}$.

When $\chi<1$, the schemes DS$_\textnormal{C,W,R}$ and AS$_\textnormal{C,W,R}$,
that use $\mathbb{G}^\textnormal{C,W,R}$,
converge after a number of iterations not exceeding $430$.
As shown in Fig.~(\ref{fig:convrate}), the number of iterations is nearly constant in the range $0\leq\chi\leq 10^{-5}$.
As expected, the accelerated schemes AS$_\textnormal{C,W,R}$ are faster than 
the direct schemes DS$_\textnormal{C,W,R}$,
with scheme AS$_\textnormal{R}$ proving the fastest.
With this scheme and when $\chi<1$, the number of iterations is at most $168$.
Again, these results are qualitatively similar with that given in~\cite{Willot14a} in the context of conductivity.

For rigidly-reinforced media ($\chi>1$), the number of iterations of schemes DS$_\textnormal{C,W,R}$
and AS$_\textnormal{C,W,R}$ follow 
the same powerlaw behaviors, with respect to $\chi$,
as that of DS or AS.
In all considered schemes, the number of iterations continuously increases with the contrast.
Differences are observed between the various
accelerated schemes AS and AS$_\textnormal{C,W,R}$, with AS$_\textnormal{R}$ the fastest.
The use of Green operators associated to the problem for the strain field, as undertaken here,
results in convergence properties that are worse in the region $\chi>1$ than when $\chi<1$.
In this respect, benefits are to be expected from the use of dual Green operators~\cite{Willot14a},
associated to the problem for the stress fields.

\subsection{Convergence rate with respect to the effective elastic moduli}
In this section, we focus on the accelerated schemes AS and AS$_\textnormal{C,W,R}$
and examine the rate of convergence of the various 
schemes with respect to the effective elastic moduli.
We consider the same Boolean microstructure as given in Sec.~(\ref{sec:rate})
but discretized on higher resolution grids of $256^3$ and $512^3$ voxels.
The volume fraction of the spheres are respectively $16.82$\% and $16.85$\%.
For simplicity, the contrast of properties take on two values $\chi=10^{-4}$ and $\chi=10^4$,
so that the spheres are quasi-porous or quasi-rigid.

\begin{figure}
\begin{center}
\includegraphics[width=9cm]{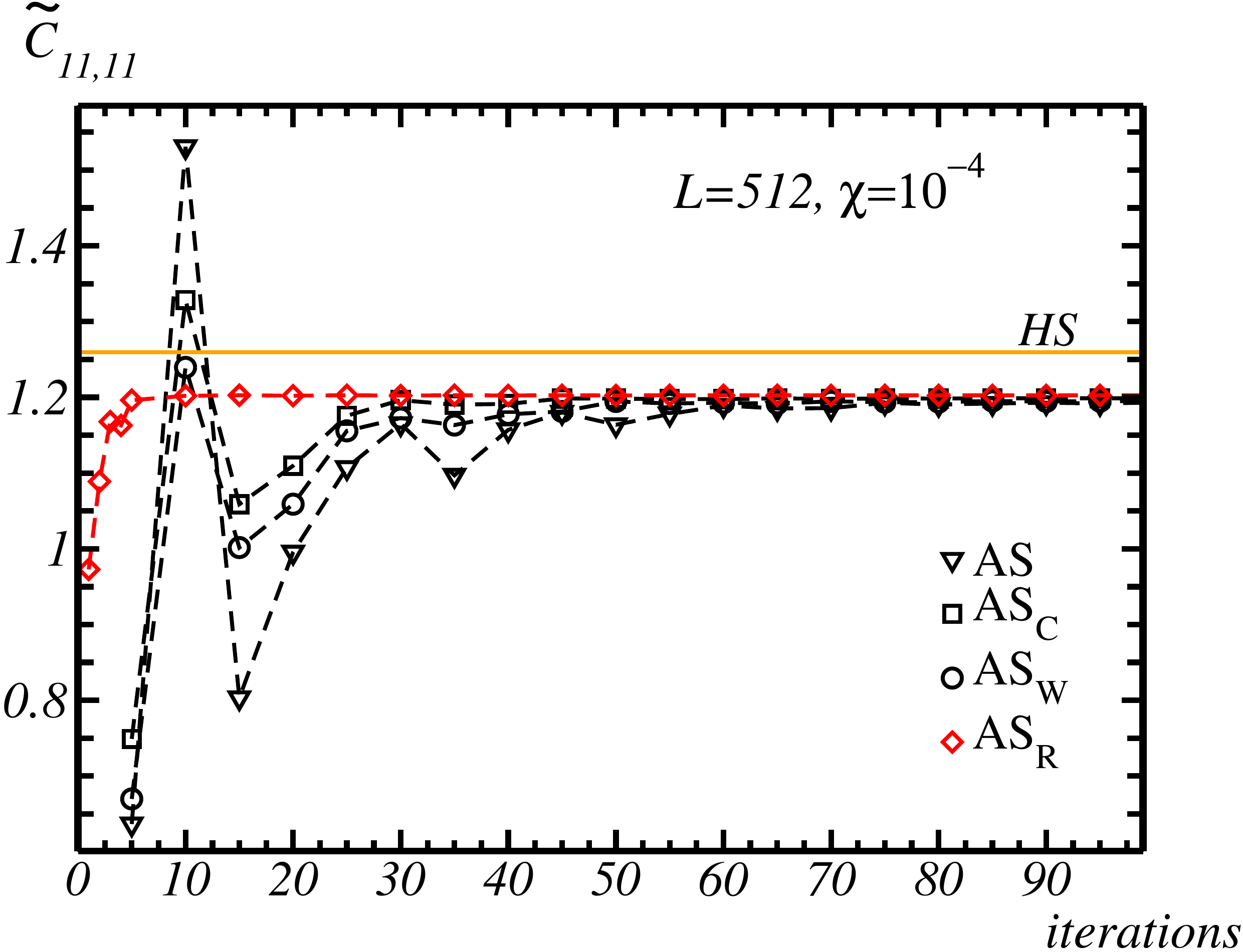} 
\caption{\label{fig:itPorous}
Estimate of the elastic modulus $\widetilde{C}_{11,11}$ as a function of the number iterations
performed, for a 3D Boolean model of quasi-porous spheres.
Black symbols: accelerated schemes AS, AS$_\textnormal{C}$, AS$_\textnormal{W}$;
red: scheme AS$_\textnormal{R}$ (orange: Hashin and Shtrikman upper bound).
}
\end{center}
\end{figure}

\begin{figure}
\begin{center}
\includegraphics[width=9cm]{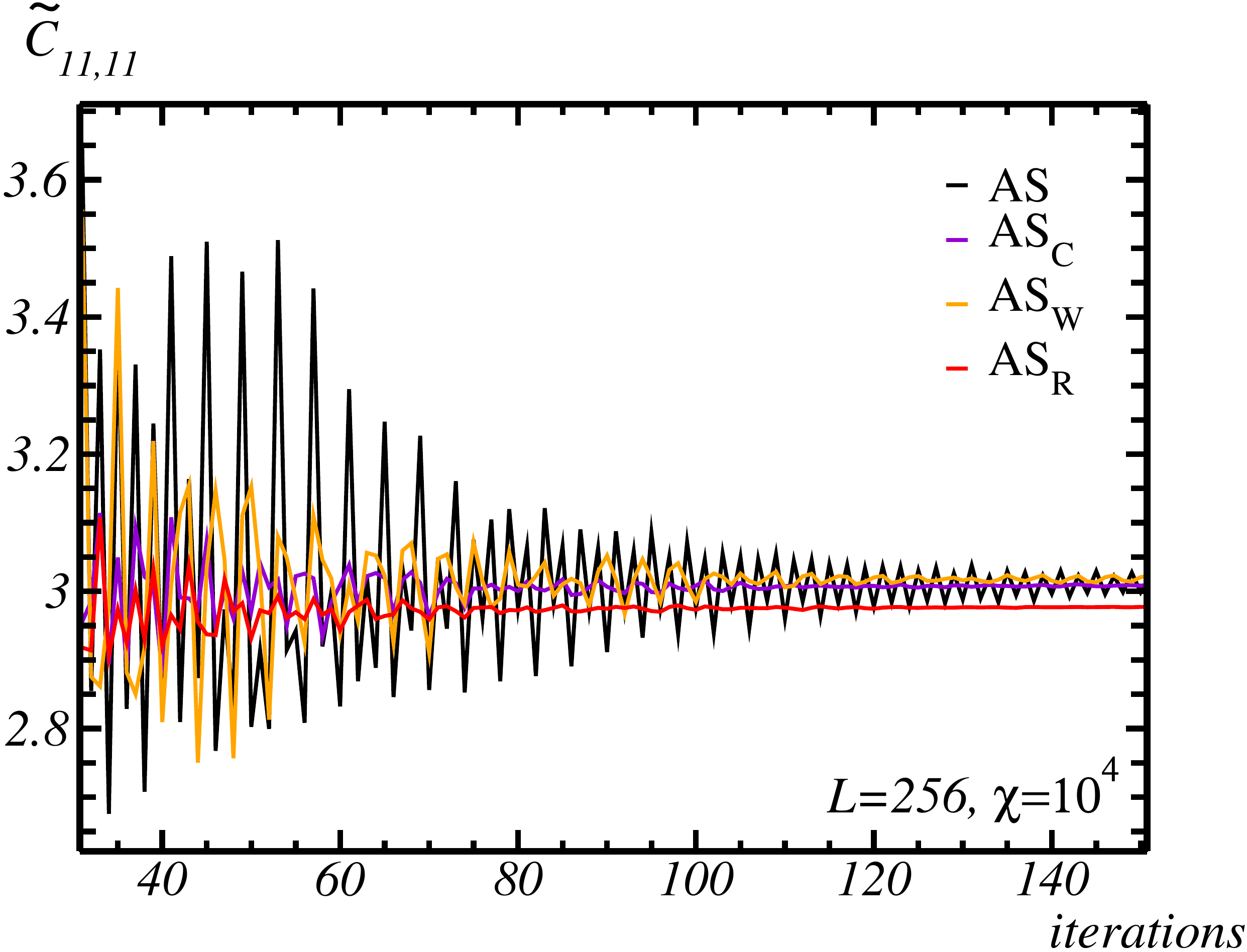} 
\caption{\label{fig:itRigid}
Estimate of the elastic modulus $\widetilde{C}_{11,11}$ as a function of the number iterations
performed, for a 3D Boolean model of quasi-rigid spheres.
Black lines: accelerated schemes AS, AS$_\textnormal{C}$, AS$_\textnormal{W}$;
red: scheme AS$_\textnormal{R}$.
}
\end{center}
\end{figure}

We perform iterations of the schemes AS and AS$_\textnormal{C,W,R}$
using the optimized reference moduli found in the previous section, on the lower resolution grid.
We apply the macroscopic strain loading:
$$
\overline{\varepsilon}_{ij}=\delta_{i1}\delta_{j1}.
$$
At each iteration and for each scheme, we 
compute the elastic modulus $\widetilde{C}_{11,11}$,
derived from the mean $\langle\sigma_{11}\rangle$ of the stress component
$\sigma_{11}$.
The convergence rate toward the elastic modulus is represented 
in Fig.~(\ref{fig:itPorous}) for quasi-porous spheres with $L=512$
and in Fig.~(\ref{fig:itRigid}) for quasi-rigid spheres with $L=256$.
In Fig.~(\ref{fig:itPorous}), for the sake of clarity, the elastic moduli are represented by symbols once every 
$5$ iterations, except for the first five iterations of the scheme AS$_\textnormal{R}$
which are all represented.
Dotted lines are guide to the eyes.
In the porous case,
much better convergence is obtained with scheme AS$_\textnormal{R}$
than with schemes AS and AS$_\textnormal{C,W}$, as shown in Fig.~(\ref{fig:itPorous}).
The estimate for $\widetilde{C}_{11,11}$
predicted by AS and AS$_\textnormal{C,W}$ present strong oscillations that are much reduced with AS$_\textnormal{R}$. 
After about $7$ iterations, the estimate given by AS$_\textnormal{R}$ is valid to a relative precision
of $10^{-2}$. 
To achieve the same precision,
more than $50$ iterations are needed for schemes AS and AS$_\textnormal{C,W}$.

The situation is notably different for quasi-rigid spheres (Fig.~\ref{fig:itRigid}).
For all schemes, a much higher number of iterations is required to determine the elastic modulus
$\widetilde{C}_{11,11}$ with a precision of $10^{-2}$.
The slower convergence rate follows that observed in Sec.~(\ref{sec:rate}),
where convergence is much poorer for $\chi>1$ than for $\chi<1$,
and where AS$_\textnormal{R}$ is less advantageous compared to the other schemes .
Nevertheless, in this case also, 
as shown in Fig.~(\ref{fig:itRigid}), smaller oscillations are observed
in the estimate for $\widetilde{C}_{11,11}$ when using AS$_\textnormal{R}$
rather than schemes AS or AS$_\textnormal{C,W}$.

\section{Conclusion}\label{sec:conclusion}
In this work, a novel discretization method has been proposed in 2D and 3D
for use in Fourier-based schemes.
The core of the proposed scheme is a simple
modification of the Green operator in the Fourier domain.
The results obtained confirm those achieved in the context of conductivity~\cite{Willot14a}. 
Compared to schemes using trigonometric polynomials
as approximation space, or to other finite-differences methods,
superior convergence rates have been observed in terms of local
stress equilibrium, but also in terms of effective properties.
More importantly, the solution for the local fields, predicted by the new discretization scheme 
is found to be more accurate than that of other methods,
especially at the vicinity of interfaces.
This property is important when applying FFT methods to
solve more complex problems like large strain deformation~\cite{Lahellec01}.
The new method also provides better estimates for the effective elastic moduli.
Furthermore, its estimates does not depend on the reference medium,
because the scheme is based on a finite-differences discretization of continuum mechanics.
Although not explored in this work, the modified Green operator 
can be used with most other FFT iterative solvers,
like the ``augmented Lagrangian''~\cite{Michel00} or with
FFT algorithms that are less sensitive to the reference~\cite{Zeman10,Brisard10},
leading to the 
same local fields.

\ack
The author thanks Carnot M.I.N.E.S for support through grant 20531.

\end{document}